\documentclass[a4paper,10pt]{amsart}

\usepackage{amsthm}                   
\usepackage{amsmath}                  
\usepackage{amsfonts}                 
\usepackage{amssymb}                  
\usepackage{mathrsfs}                 
\usepackage{xy}
\usepackage{graphicx,subfigure}
\xyoption{all}

\usepackage{epstopdf}

\theoremstyle{definition}
\newtheorem{definition}{Definition}[section]

\newtheorem{example}{Example}

\theoremstyle{plain}
\newtheorem{theorem}[definition]{Theorem}
\newtheorem{lemma}[definition]{Lemma}
\newtheorem{corollary}[definition]{Corollary}
\newtheorem{proposition}[definition]{Proposition}
\newtheorem*{theorem*}{Theorem}

\numberwithin{equation}{section}

\newcommand{\Z}{\ensuremath{\mathbb{Z}}}     

\newcommand{\bomega}{{\boldsymbol {\omega}}}
\newcommand{\bgamma}{{\boldsymbol {\gamma}}}
\newcommand{\bdeta}{{\boldsymbol {\eta}}}

\newcommand{\SH}{{\mathcal H}}

\newcommand{\e}{{\mathbf e}}
\newcommand{\ignore}[1]{}

\def\H{\mathcal H}

\def\cal{\mathcal }
\def\R{\mathbb R}

\def\Z{\mathbb Z}

\def\A{\mathcal A}
\def\G{\mathcal G}
\begin{document}
\author{Shu-qin Zhang}
\address{Chair of Mathematics and Statistics, University of Leoben, Franz-Josef-Strasse 18, A-8700 Leoben, Austria}
\email{zhangsq\_ccnu@sina.com}
\title[Optimal parametrizations of a class of self-affine sets]{Optimal parametrizations of a class of self-affine sets}

\date{\today}
\thanks{Supported by project I1136 and by the doctoral program W1230 granted by the Austrian Science Fund (FWF)}

\maketitle

\begin{abstract}
In this paper, we study optimal parametrizations of the invariant sets of a single matrix graph IFS which is a generalization of \cite{RaoZhangS16}.We show that the invariant sets of a linear single matrix GIFS which has primitive associated matrix and  satisfies the open set condition admit optimal parametrizations. This result is the basis of further study of space-filling curves of self-affine sets.
\end{abstract}

\section{Introduction }
The topic of space-filling curves has a very long history. Recently,
Rao and Zhang \cite{RaoZhangS16} as well as Dai, Rao, and Zhang \cite{DaiRaoZhang15} found a systematic method
to construct space-filling curves for connected self-similar sets satisfying the open set condition.
This method generalizes almost all known results in this field.

The self-affine sets have more complex structures than self-similar due to the different contraction ratios in different directions. There are almost no systematic works on the space-filling curves of self-affine sets except some examples provided by
Dekking \cite{Dekking82b},  Sirvent's study under some special conditions \cite{Sirvent2008,Sirvent2012}, boundary parametrizations of self-affine tiles by Akiyama and Loridant \cite{AkiyamaLoridant10, AkiyamaLoridant11}, and boundary parametrizations of a class of cubic Rauzy fractals by Loridant \cite{Loridant2016}.
The purpose of the present paper is to carry out studies in this direction. 
In particular, we study the space-filling curves (or optimal parameterizations) of single-matrix graph-directed systems.

\subsection{Single-matrix GIFS}
Let $(\A,\Gamma)$ be a directed graph with vertex set $\A$ and edge set $\Gamma$.
Let
$$\mathcal{G}=\{g_e: \R^d\rightarrow \R^d;~e\in \Gamma\}$$ be a collection of contractions.
The triple $(V,\Gamma,\mathcal{G})$, or simply $\mathcal{G}$, is called \emph{a graph-directed iterated function system }(GIFS). Usually, we set $\A=\{1,2,\dots,N\}$.
Let $\Gamma_{ij}$ be the set of edges from vertex $i$ to $j$, then there exist unique non-empty compact sets $\{E_i\}_{i=1}^N$ satisfying
\begin{equation}\label{GIFS1}
E_i=\bigcup_{j=1}^N \bigcup_{e\in\Gamma_{ij}} g_e(E_j),\quad\quad 1\leq i \leq N.
\end{equation}
The family $\{E_i\}_{i=1}^N$ is called the family of \emph{invariant sets} of the GIFS (cf. \cite{MauldinWilliams88}).
When the graph has only one vertex with self-edges, then the GIFS will degenerate into an iterated function system (IFS).

We say that the system $\mathcal{G}$ satisfies the \emph{open set condition} (OSC) if there exist non-empty open sets $U_1,\dots, U_N$ such that 
$$U_i\subset \bigcup_{j=1}^N \bigcup_{e\in\Gamma_{ij}} g_e(U_j),\quad\quad 1\leq i \leq N,
$$
and the right hand sets union are disjoint.

Denote $M=(m_{ij})_{1\leq i,j\leq N}$ the \emph{associated matrix} of the directed graph $(\A,\Gamma)$, that is, $m_{ij}=\sharp\Gamma_{ji}$ counts the number of edges from $j$ to $i$. We say a directed graph $({\mathcal A}, \Gamma)$ is \emph{primitive}, if the associated matrix is primitive, i.e., $M^n$ is a positive matrix for some $n$.
(See \cite{MauldinWilliams88}, \cite{Falconer90}.)  In this paper, we always assume that $({\mathcal A}, \Gamma)$ is primitive.



According to Luo and Yang \cite{LuoYang10}, ${\mathcal G}$ is called a \emph{single-matrix GIFS} if there is a $d\times d$ expanding matrix $A$ such that all functions related to $e\in \Gamma$ have the form
\begin{equation}\label{GIFS2}
g_e(x)=A^{-1}(x+d_e),
\end{equation}
where $d_e \in \mathbb{R}^d$.





\subsection{Optimal parametrization}
Let $E\subset \R^d$ be a compact set and $\SH^s(E)$ denote the Hausdorff measure with respect to Euclidean norm of $E$. Basically, if
$\psi:[0,1]\to E$ is a continuous onto mapping, then $\psi$ is a \emph{parametrization} of $E$. 
If $E$ is a self-similar set satisfying the open set condition, then
$0<\SH^s(E)<\infty$, where $s$ is the Hausdorff dimension of $E$.
In this case, we may expect that $E$ has a better parametrization.
The following concept is first given by Dai and Wang \cite{DaiWang2010}:

\begin{definition}\label{optimal} {\rm
 A surjective mapping   $\psi:[0,1]\rightarrow E$ is called an
\emph{optimal parametrization} of $E$ if the following conditions are fulfilled.
\begin{enumerate}
  \item[($i$)] 
   $\psi$ is a \emph{measure isomorphism} between $([0,1],$$ {\cal B}([0,1]),$$ {\mathcal L})$ and
   $(E, {\cal B}(E),\SH^s)$, that is,  there exist $E'\subset E$ and $I'\subset [0,1]$ with full measure such that
   $\psi:~I'\to E'$ is a bijection and it is  measure-preserving in the sense that
  $$
  {\mathcal H}^s(\psi(F))=c{\mathcal L}(F) \text{ and } {\mathcal L}(\psi^{-1}(B))=c^{-1}{\mathcal H}^s(B),
  $$
  for any Borel set $F\subset [0,1]$ and any Borel set $B\subset E$, where $c={\mathcal H}^s(E)$. (See for instance, Walters \cite{Walters1982}.)
  \item[($ii$)] $\psi$ is \emph{$1/s$-H\"older continuous}, that is, there is a constant $c'>0$ such that
  $$
  \|\psi(x)-\psi(y)\|\leq c' \|x-y\|^{\frac{1}{s}} \  \text{ for all } x,y\in [0,1].
  $$
 We call $1/s$ the \emph{H\"older exponent}.
  \end{enumerate}
  
}
\end{definition}

For a self-affine set $K$, the Hausdorff measure may be $0$ or $\infty$, and hence
we cannot require an optimal parametrization satisfying (i) of the above. Also, the $1/s$-H\"older continuity may fail.
So we are forced to define the optimal parametrization in some other way.

To this end, we choose a pseudo-norm $\|\cdot\|_\omega$ instead of the Euclidean norm on $K$.
This pseudo-norm is first introduced by
Lemari{\'e}-Rieusset \cite{LemarieRieusset94} to deal with problems in the theory of wavelets.
   He and Lau \cite{HeLau08} developed the Hausdorff dimension (denoted by $\text{dim}_\omega$) and Hausdorff measure (denoted by $\H_{\omega}^s$) with respect to pseudo-norm (see Section \ref{Preli} for details). The advantage of the pseudo-norm is that we can regard the expanding matrix $A$ as a `similitude'.
   By replacing the norm, dimension and measure by their counterpart w.r.t.the pseudo-norm, we can define an optimal parametrization similar to Definition \ref{optimal}; details will be given in Section \ref{Preli}.

\subsection{Linear GIFS and our main result}
We equip a GIFS with an order structure (and
call it an \emph{ordered GIFS}), which induces a lexicographical order of the associated symbolic
space. An ordered GIFS is called a \emph{linear GIFS} if every two consecutive cylinders have nonempty intersections (see section \ref{Preli} for precise definitions).

Rao and Zhang \cite{RaoZhangS16} proved that as soon as we find a linear GIFS structure of a self-similar set, then
a space-filling curve can be constructed accordingly. Dai, Rao, and Zhang \cite{DaiRaoZhang15}
develop a very general method to explore linear GIFS structures of a given self-similar set.
In this paper, we only deal with the first problem for self-affine sets. Our main result is

\begin{theorem} \label{Main1}
Let $(\mathcal{A}, \Gamma, \mathcal{G},\prec)$ be  a linear single-matrix graph-directed IFS with expanding matrix $A$ satisfying the open set condition and assume that the associated matrix $M$ of the graph is primitive. Then there exists a parametrization $\psi_j$ of the invariant $E_j$ for all $j\in \mathcal{A}$ such that
\begin{enumerate}
  \item[($i$)] 
   $\psi_j$ is a \emph{measure isomorphism} between $$([0,1], {\cal B}([0,1]), {\mathcal L}) \text{ and } (E_j, {\cal B}(E_j),\SH_{\omega}^s).$$
  \item[($ii$)] There is a constant $c>0$ such that
  $$
  \|\psi_j(x)-\psi_j(y)\|_\omega\leq c\|x-y\|^{\frac{1}{\alpha}} \  \text{ for all } x,y\in [0,1],
  $$
where $\alpha=\dim_{\omega} E_j$. 
 \end{enumerate} 

\end{theorem}

According to the relation between Euclidean norm and pseudo-norm, we have the following result for the H\"older continuity of the parametrization $\psi_j$ obtained by the above theorem.
\begin{corollary}\label{intro-cor}
 Let $\lambda_M$ be the maximal eigenvalue of $M$. Let $\lambda_{max}$ and $\lambda_{min}$ be the maximal modulus and minimal modulus of the eigenvalues of $A$, respectively. 
For any $0<\epsilon<\lambda_{min}-1$,

We will close this section by three examples, where the substitution rules induced linear GIFS are given directly. We shall discuss
in subsequent paper on how to obtain substitution rules, and once we have a substitution rule, the GIFS is constructed automatically. In the following examples, we calculate the value of H\"older exponent and it will be interesting to compare this with the value obtained from Corollary \ref{intro-cor}.

\begin{itemize}
\item  $\psi_j$ is $\frac{\ln (\lambda_{\max}+\epsilon)}{\ln \lambda_M}$-H\"older continuous if $\|\psi_j(x)-\psi_j(y)\|\geq 1$;

\item $\psi_j$ is $\frac{\ln (\lambda_{\min}-\epsilon)}{\ln \lambda_M}$-H\"older continuous if $\|\psi_j(x)-\psi_j(y)\|\leq 1$.
\end{itemize}
 The matrices $M$ and $A$ have the same meaning as in the above theorem.
\end{corollary}
\begin{figure}[htbp]
\subfigure[\text{Unit square $Q=\cup_{i=1}^6S_i(Q)$.}]{\includegraphics[width=0.33 \textwidth]{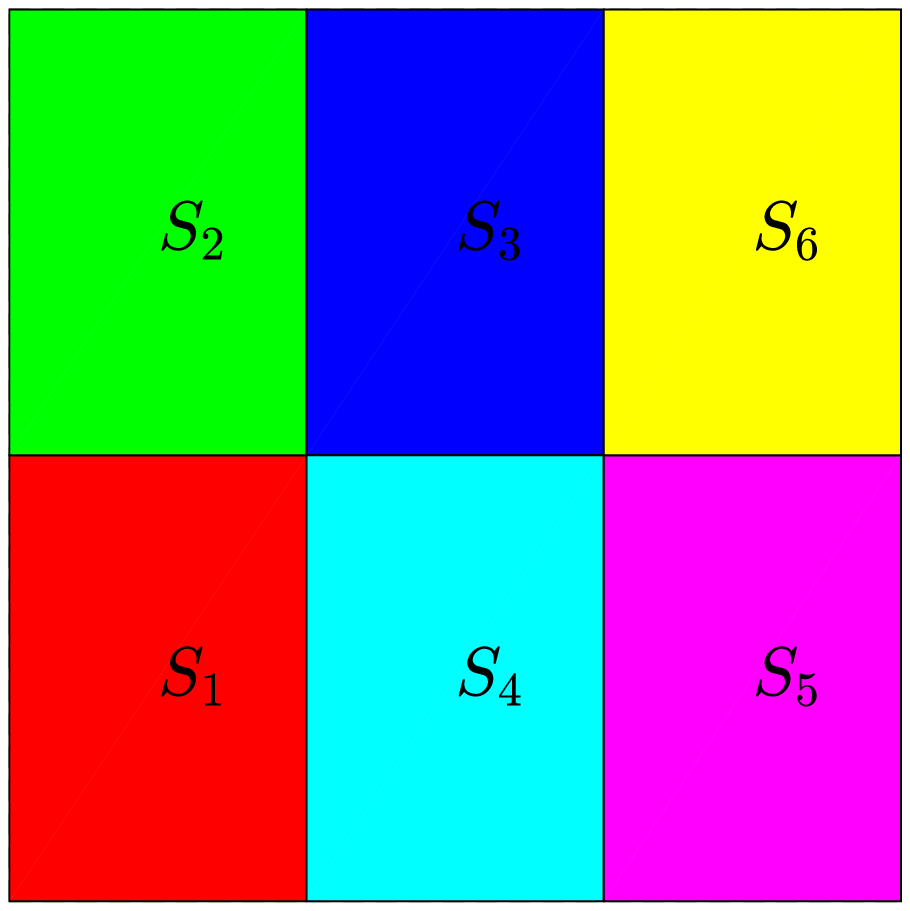}}
\subfigure[\text{Initial cycle.}]{\includegraphics[width=0.33  \textwidth]{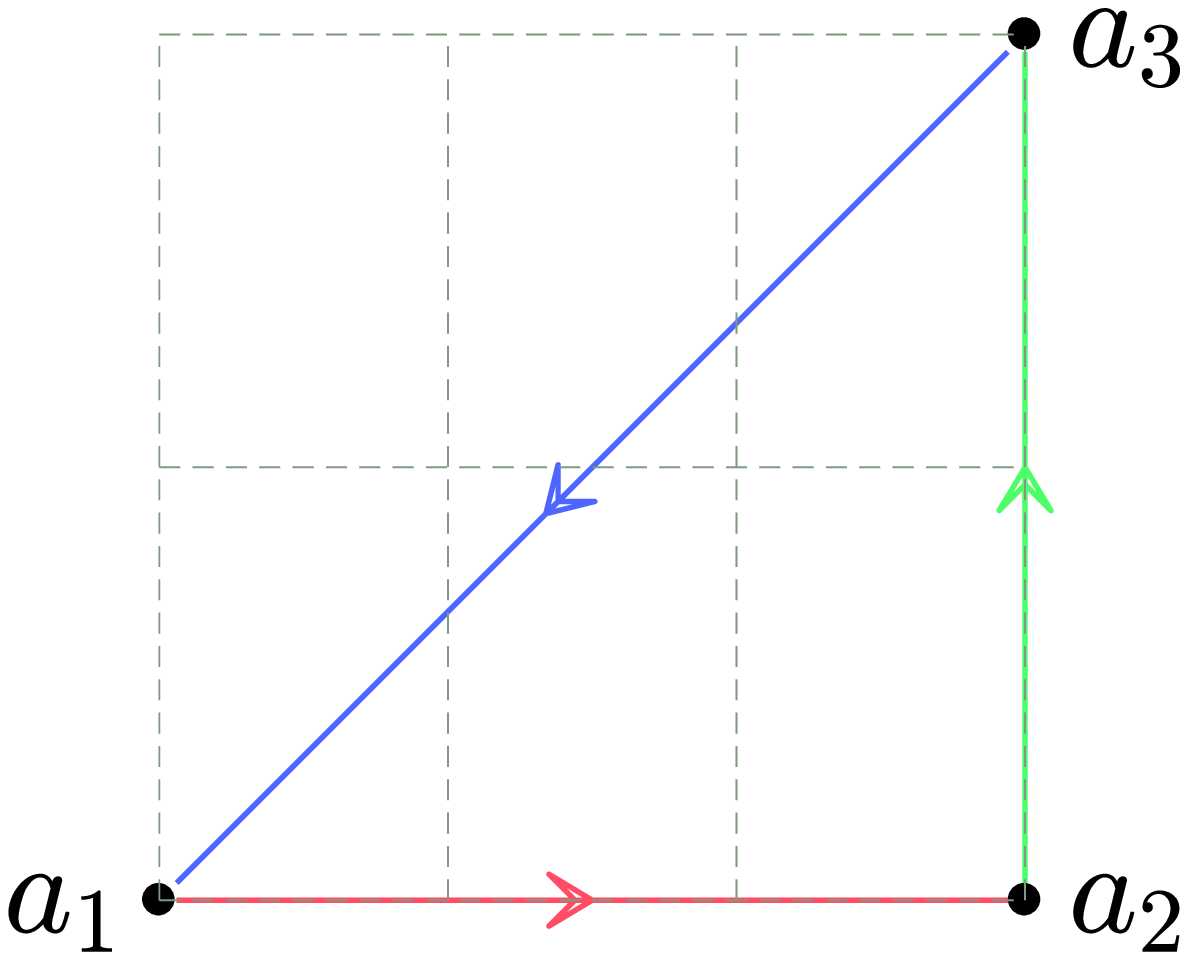}}\subfigure[\text{First generation.}]{\includegraphics[width=0.33  \textwidth]{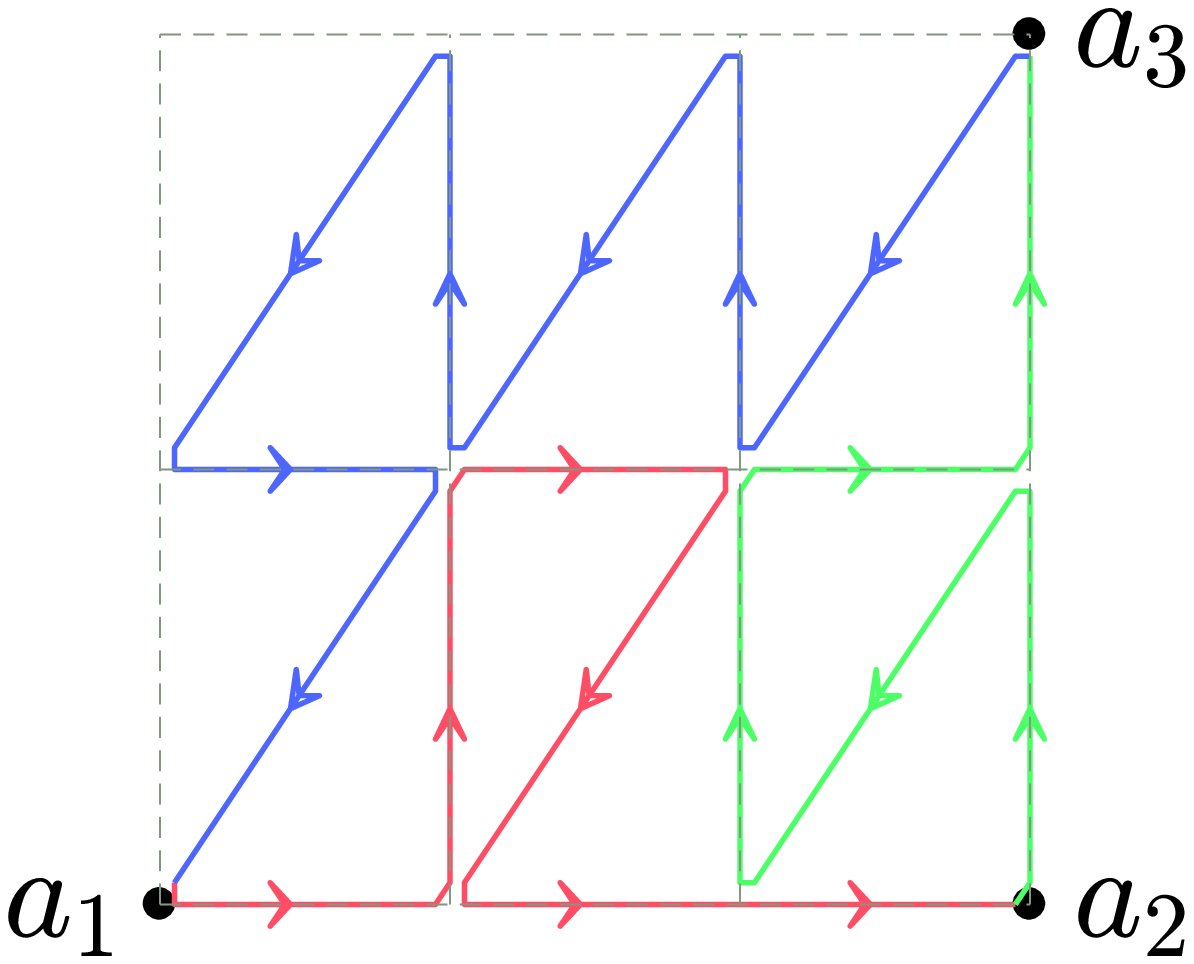}}\\
\subfigure[\text{Second generation.}]{\includegraphics[width=0.33 \textwidth]{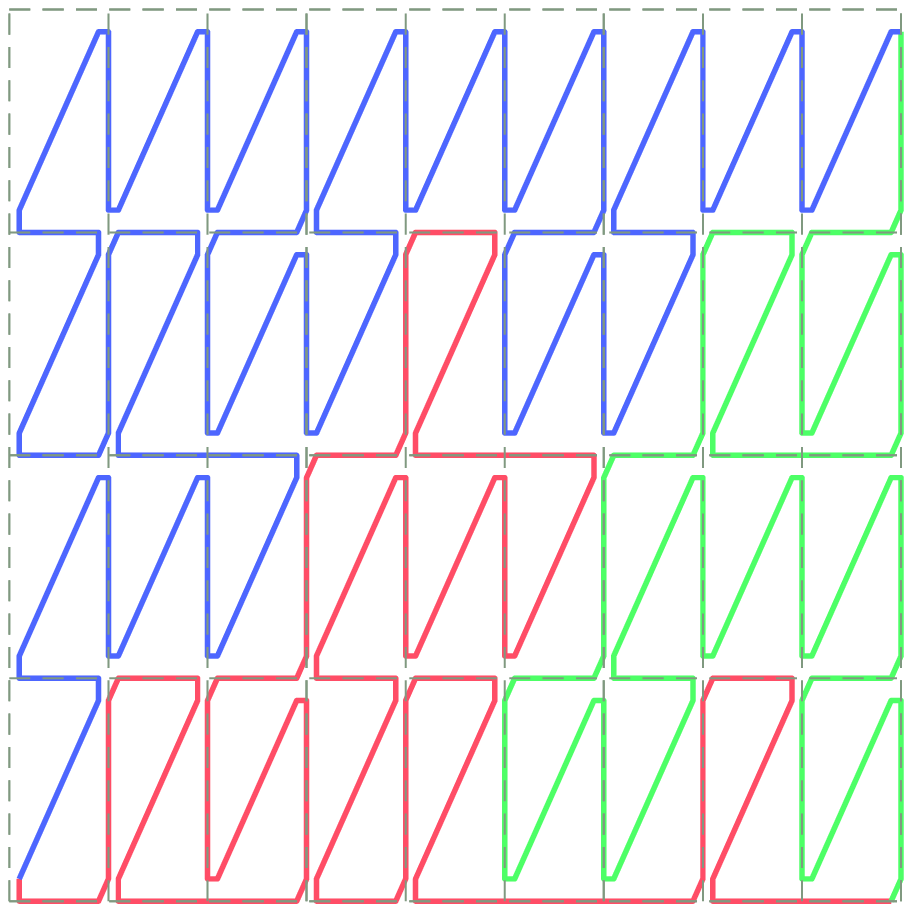}}\subfigure[\text{Third generation.}]{\includegraphics[width=0.33 \textwidth]{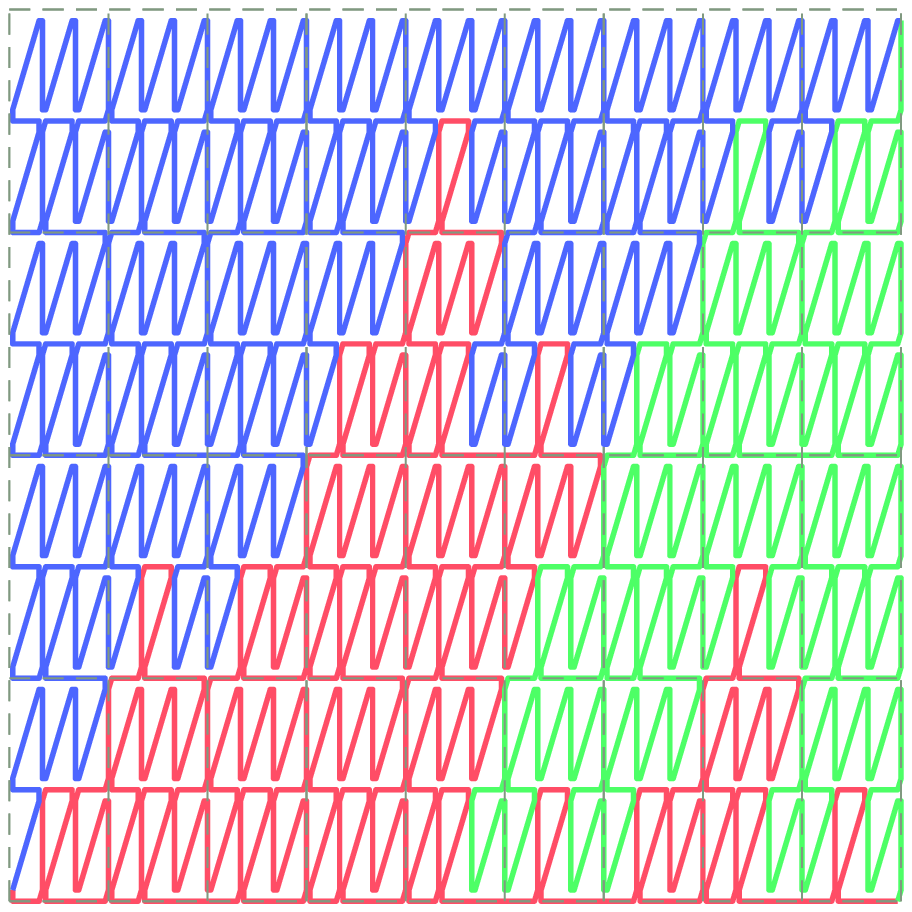}}\subfigure[\text{ Parametrized square.}]{\includegraphics[width=0.33  \textwidth]{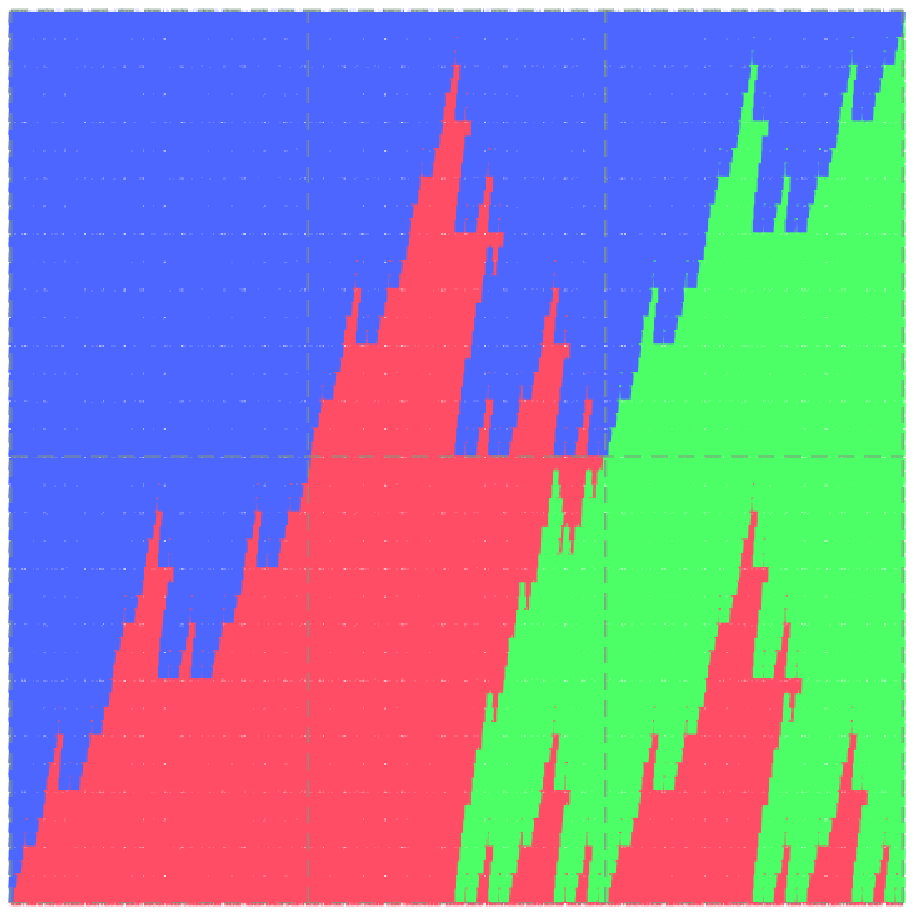}}
\caption{The space-filling curve of unit square given by IFS $\{S_i\}_{i=1}^6$.}\label{Ex3_fig1}
\end{figure}

In the following examples, we always use ${\bf e}_1, {\bf e}_2$ to denote the standard basis of $\mathbb{R}^2$.  Denote the maximal eigenvalue of $M$ by $\lambda_M$. And denote by the H\"older exponent with respect to Euclidean norm H\"older$_E$.

\begin{example}\label{Ex3}
\begin{figure}[htbp]
\includegraphics[height=7.5 cm]{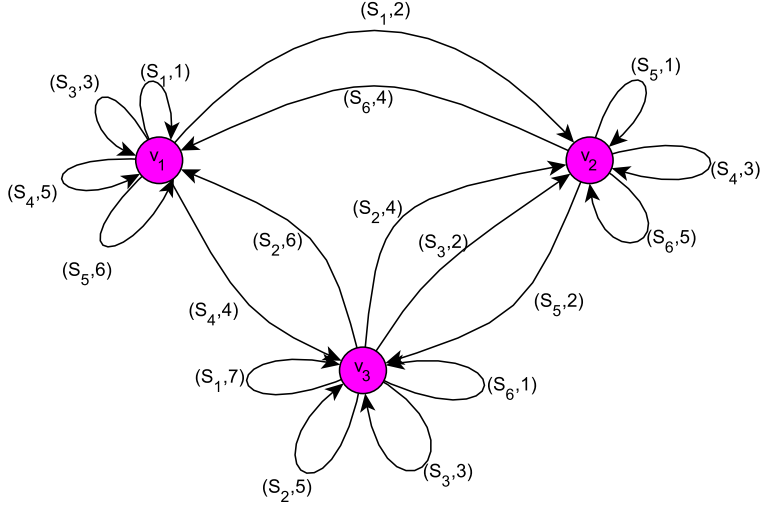}
\caption{The directed graph $G$ with vertex set $\{v_1, v_2, v_3\}$ and labelled edges obtained by the substitution rule \eqref{Ex3_subrule}. The label $(S,i)$ with $S\in \{S_i\}_{i=1}^6$ and $i\in \mathbb{N}$ means that $S$ is the contraction for the $i$-th edge starting from this vertex. The graph here determins a GIFS.}\label{GIFS_graph}
\end{figure}
{\bf A unit square.}  Let $Q$ be the unit square generated by the IFS $S_i(x)=A^{-1}(x+d_i),~~d_i\in \mathcal{D}$, where $\mathcal{D}=\{0,{\bf e}_2, {\bf e}_1+{\bf e}_2, {\bf e}_1, 2{\bf e}_1,2{\bf e}_1+{\bf e}_2\}$, and the expanding matrix is $A=\left(\begin{matrix}
3&0\\
0&2\\
\end{matrix}\right)$.
Let $a_1,a_2,a_3$ be the three vertice of the unit square, see Figure \ref{Ex3_fig1}(b). Denote by $v_i$ the edge from $a_i$ to $a_{i+1}, i=1, 2, 3$ (assume $a_4=a_1$). By the similar idea as Dai, Rao, and Zhang \cite{DaiRaoZhang15}, we can construct the following substitution rule which means to replace an edge $v_i$ by a trail which shares the same starting point and ending point with $v_i$.
\begin{equation}\label{Ex3_subrule}
\begin{split}
v_1 \longrightarrow\vphantom{a}&S_1(v_1)+ S_1(v_2)+S_3(v_1)+S_4(v_3)+S_4(v_1)+S_5(v_1),
\\
v_2 \longrightarrow\vphantom{a}&S_5(v_2)+ S_5(v_3)+S_4(v_2)+S_6(v_1)+S_6(v_2),
\\
v_3 \longrightarrow\vphantom{a}&S_6(v_3)+ S_3(v_2)+S_3(v_3)+S_2(v_2)+S_2(v_3)+S_2(v_1)+S_1(v_3),
\end{split}
\end{equation}
where we use the symbol $'+'$ to connect the consecutive edges or sub-trails. 
It is shown in \cite{DaiRaoZhang15} that the above substitution rule induces a GIFS (see Figure \ref{GIFS_graph}) which is linear. We will show the set equation form of the GIFS in Section \ref{Sec_Example}.
The associated matrix of the substitution which is defined as the associated matrix of the directed graph $G$ obtained by the substitution is $$M=\left(\begin{matrix}
4&1&1\\
1&3&2\\
1&1&4\\
\end{matrix}\right),\quad\lambda_M=6,\quad\text{H\"older}_E=\log_6 2.$$ 
Compared with the unit square parametrized using the method as Hilbert or Peano which have the H\"older exponent $\frac{1}{2}$, the parametrization obtained here doesn't have a better smoothness.
\end{example}

\begin{figure}[htbp]
\includegraphics[width=0.38 \textwidth]{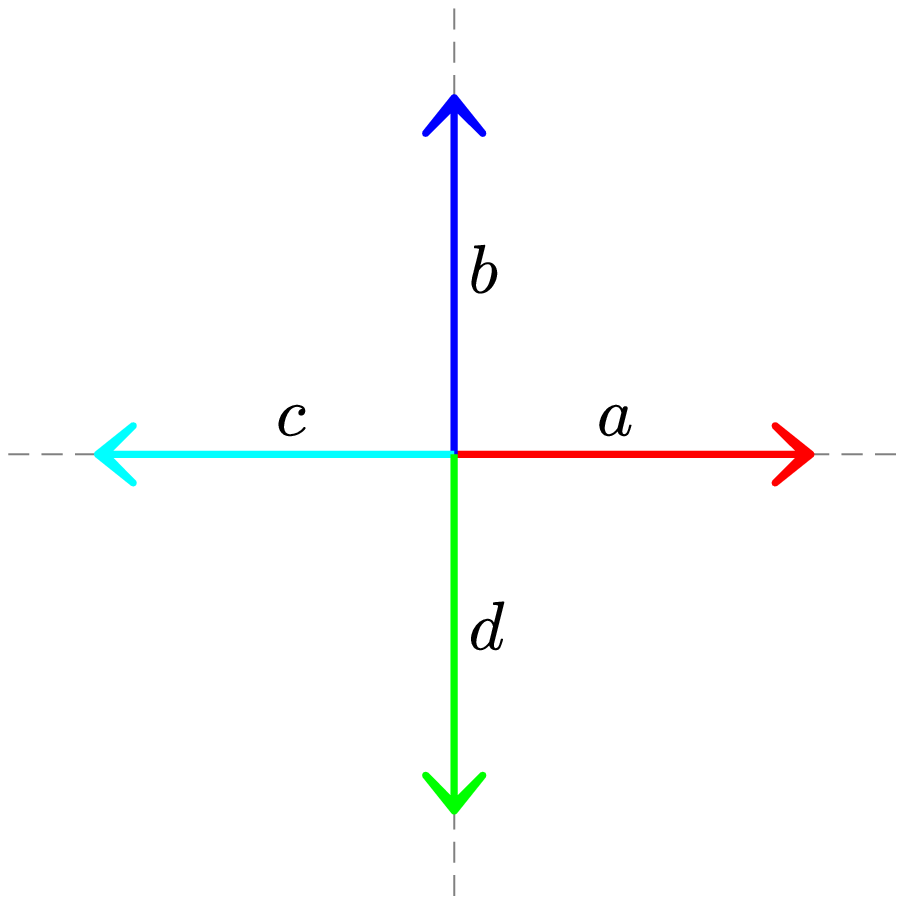}
\includegraphics[width=0.3 \textwidth]{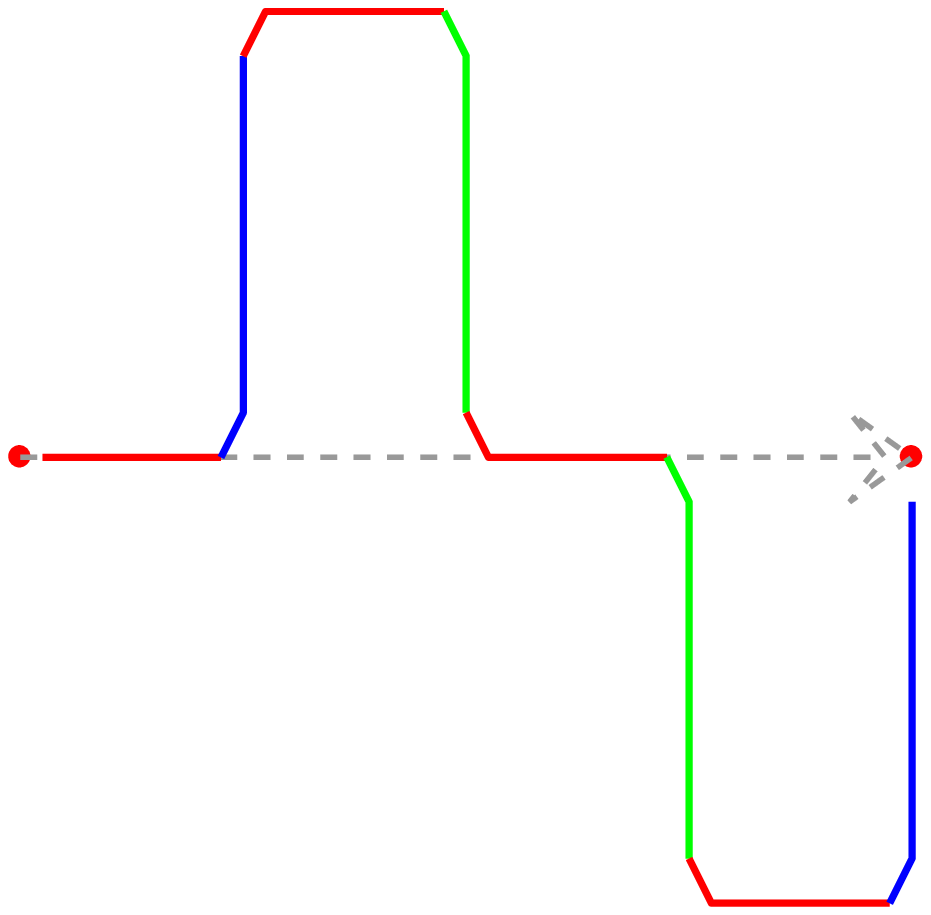}\includegraphics[width=0.3 \textwidth]{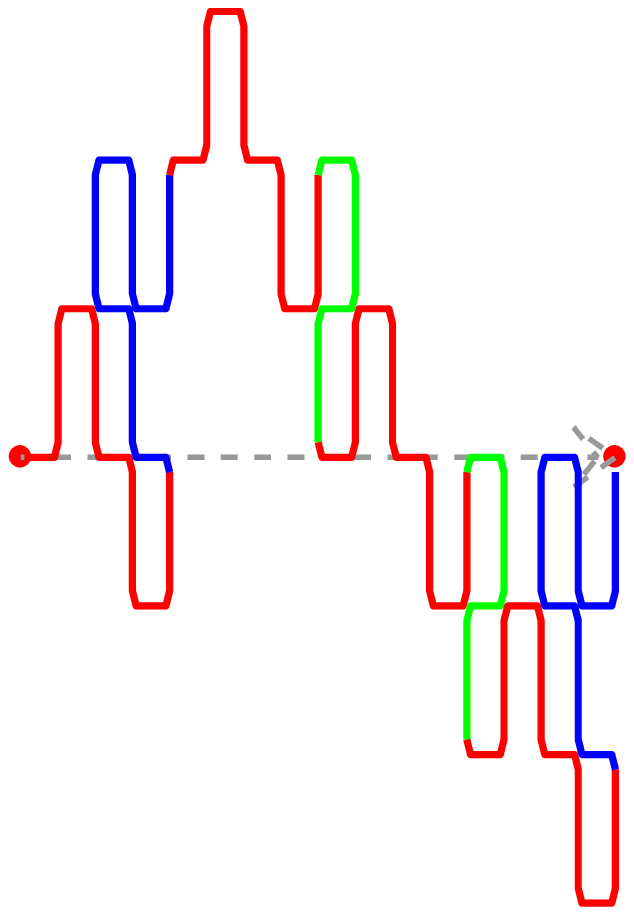}\\
\includegraphics[width=0.4 \textwidth]{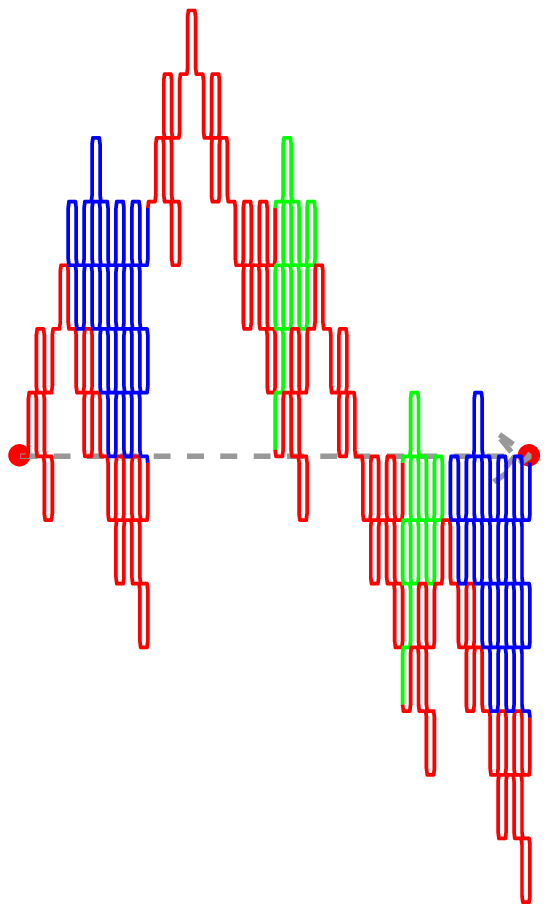}
\includegraphics[width=0.4 \textwidth]{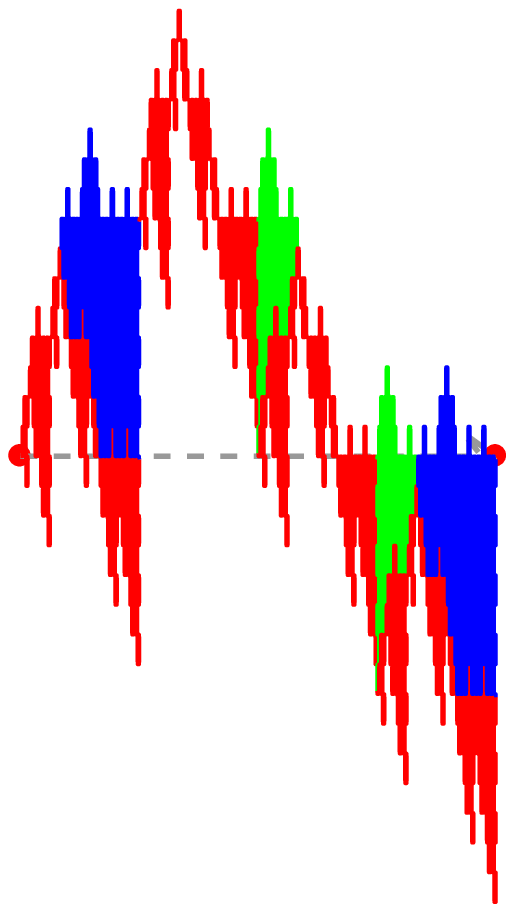}
\caption{The figures show the first four approximation of $E_1$.}\label{Ex1_Fig2}
\end{figure}

\begin{example}\label{Ex1}
{\bf Dekking's plane filling curve \cite{Dekking82b}}. It is induced by the following substitution:
$$\begin{matrix}
\sigma: & a\mapsto abadadab; \quad b\mapsto cbcbadab; \quad c\mapsto cbcbcdadcbcd;\quad d\mapsto adcd.
\end{matrix}
$$
The set equation form of this substitution can be found in Section \ref{Sec_Example} and we obtain a linear GIFS as well. In this example, the expanding matrix is $A=\left(\begin{matrix}
4&0\\
0&2\\
\end{matrix}\right)$ and associated matrix of the substitution is
$$M=\left(\begin{matrix}
4&2&1&1\\
2&3&3&0\\
0&2&5&1\\
2&1&3&2\\
\end{matrix}\right),  \quad\lambda_M=8,\quad\text{H\"older}_E=\frac{1}{3}.$$
Then we can check that H\"older$_E$ is between the two H\"older exponents obtained by Corollary \ref{intro-cor}.
\end{example}


\begin{example}\label{Ex2}
\begin{figure}[h]
\includegraphics[width=0.35 \textwidth]{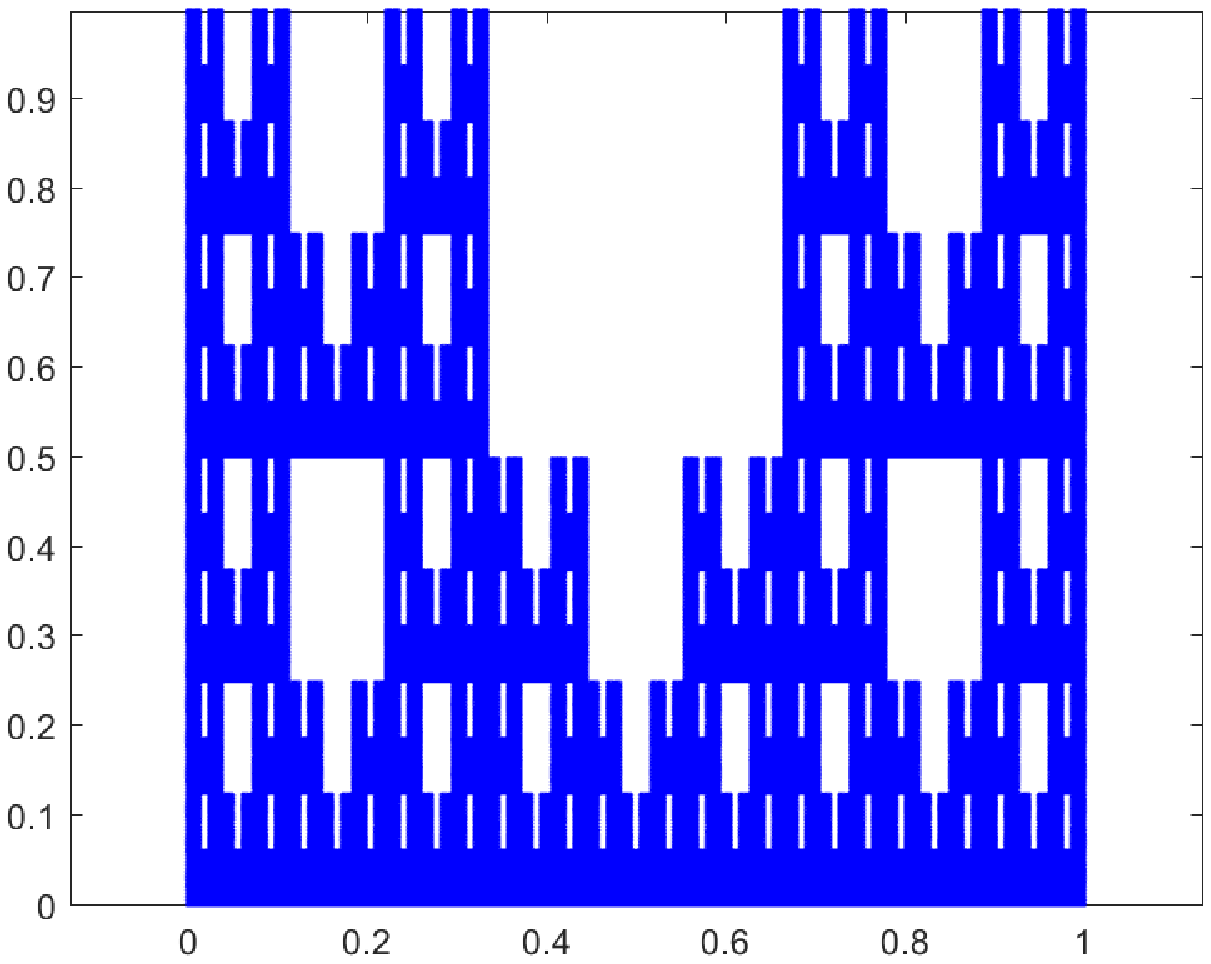}\quad\quad \includegraphics[width=0.35 \textwidth]{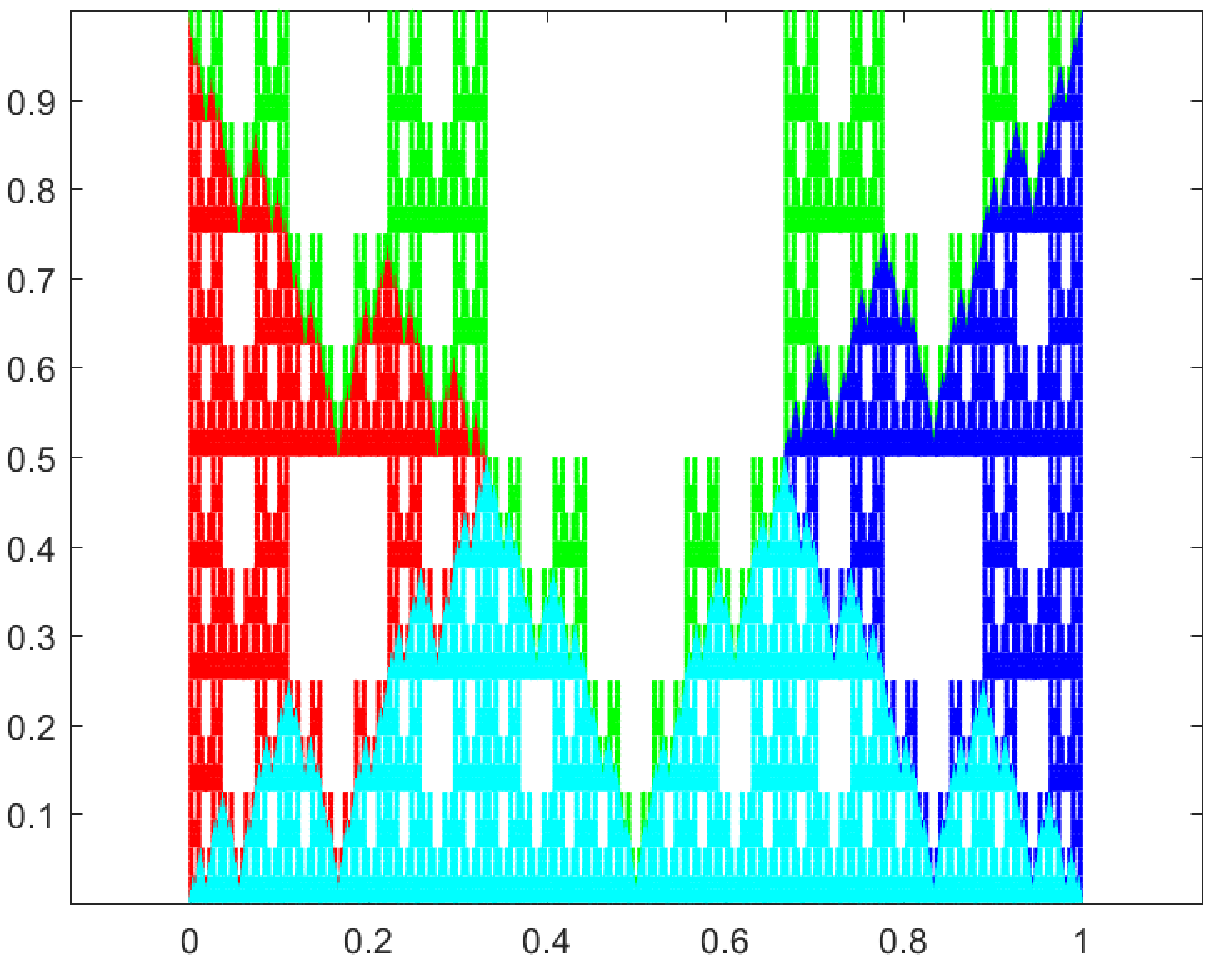}
\caption{The left is Mcmullen set $T$ and the right is $\cup_{i=1}^5 E_i$.}\label{Ex2_fig1}
\end{figure}
\begin{figure}[h]
\includegraphics[width=0.35 \textwidth]{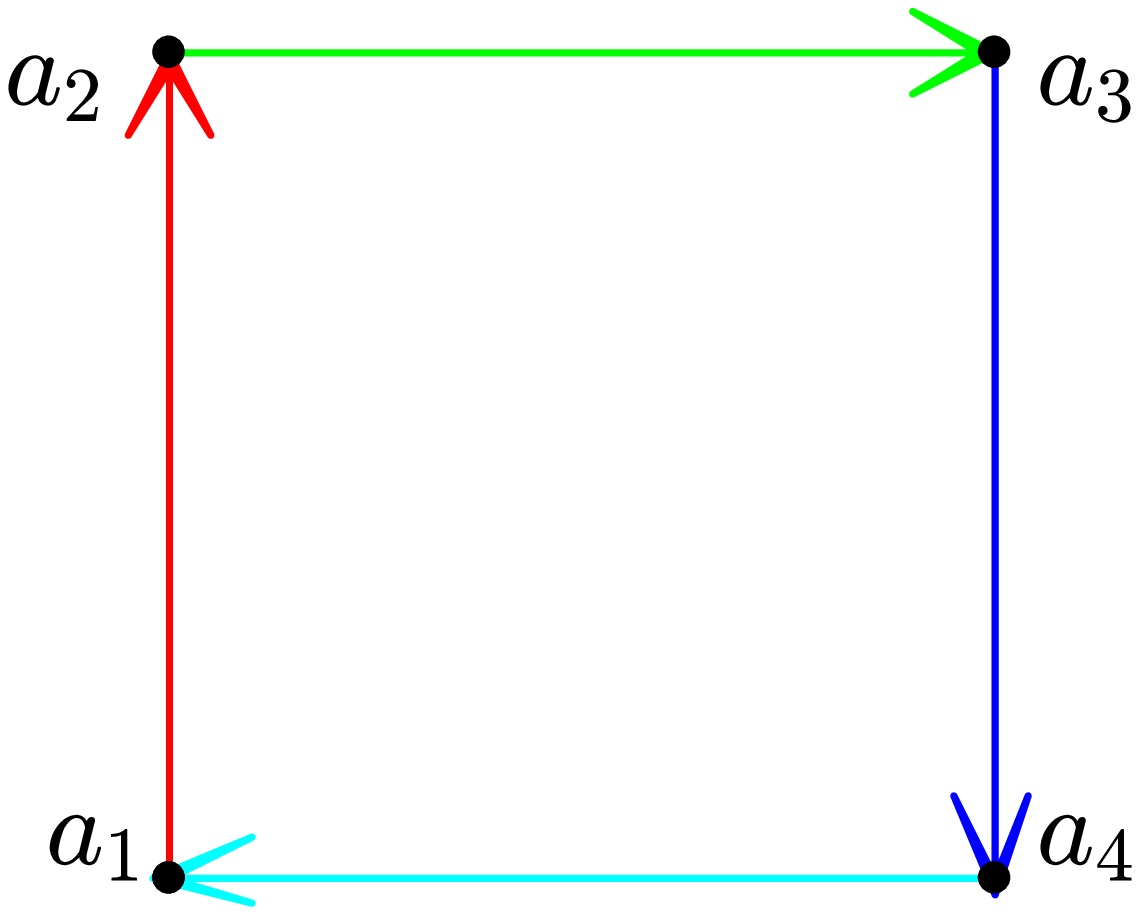}\includegraphics[width=0.35 \textwidth]{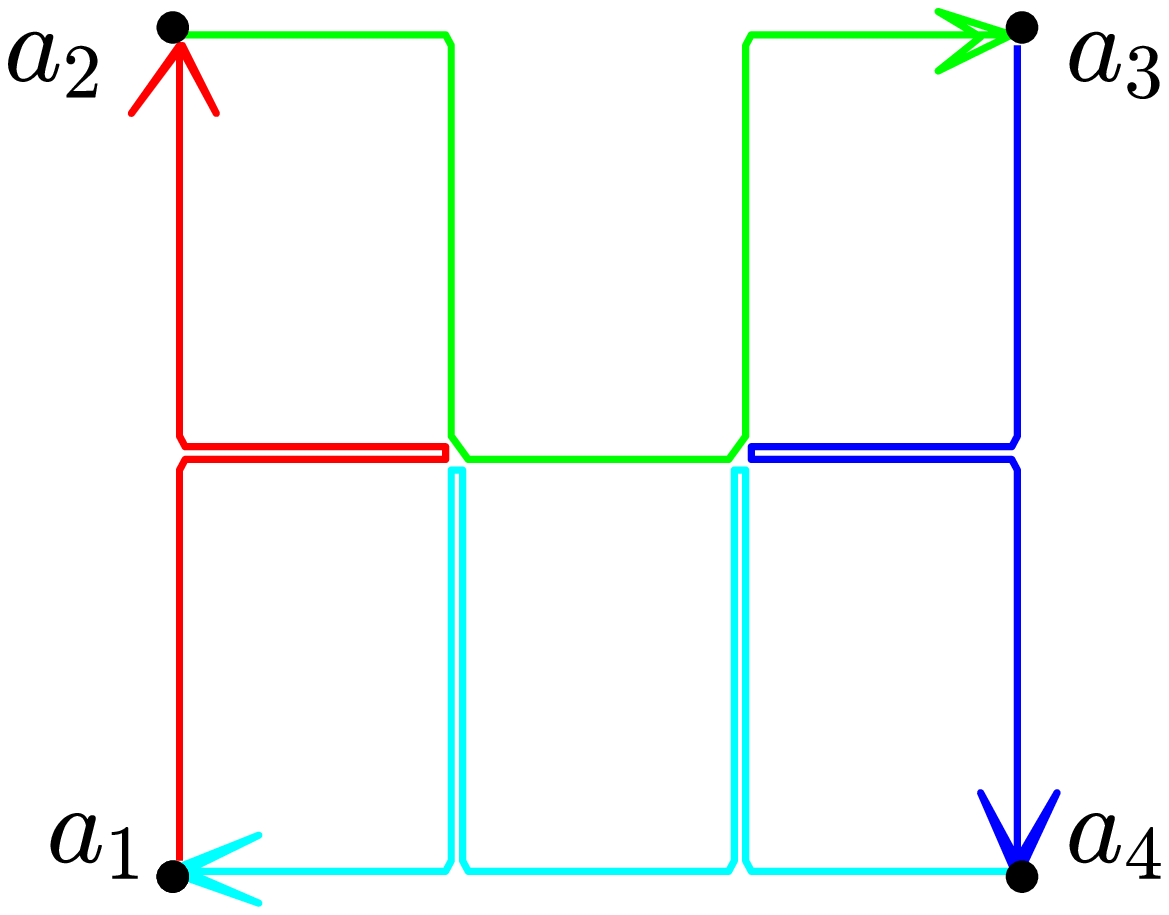}
\caption{Substitution rule of the Mcmullen set }\label{Ex2_fig2}
\end{figure}
{\bf A McMullen set}.  We consider the McMullen set depicted in Figure~ \ref{Ex2_fig1}, left. Denote the four vertices of the unit square by $a_1, a_2, a_3, a_4$.
Denote
$v_i=\overrightarrow{a_ia_{i+1}}, i=1,2,3,4~ (\text{ assume } a_5=a_1),$
and we get the following substitution rule (which is obtained from replacing the edge from $a_i$ to $a_{i+1}$ (left, Figure \ref{Ex2_fig2}) by the trail from $a_i$ to $a_{i+1}$ (right, Figure \ref{Ex2_fig2})).
\begin{flalign*}
 v_1 \longrightarrow\vphantom{a}&S_1(v_1)+ S_1(v_2)+S_4(v_4)+S_4(v_1),
\\
v_2 \longrightarrow\vphantom{a}&S_4(v_2)+ S_4(v_3)+S_2(v_2)+S_5(v_1)+S_5(v_2),
\\
v_3 \longrightarrow\vphantom{a}&S_5(v_3)+ S_5(v_4)+S_3(v_2)+S_3(v_3),
\\
v_4\longrightarrow\vphantom{a}&S_3(v_4)+S_3(v_1)+S_2(v_3)+S_2(v_4)+S_2(v_1)+S_1(v_3)+S_1(v_4).
\end{flalign*}
In this example, the expanding matrix is $A=\left(\begin{matrix}
3 &0\\
0 &2\\
\end{matrix}\right)$ and the associated matrix of the substitution is $$M=\left(\begin{matrix}
2&1&0&2\\
1&3&1&0\\
0&1&2&2\\
1&0&1&3\\
\end{matrix}\right), \quad\lambda_M=5,  \quad \text{H\"older}_E=\log_5 2.$$
Figure
\ref{Ex2_fig3} shows the visualization of the filling curve of the Mcmullen set $T$. To give a self-avoiding visualization, we round off the corners of the approximating curves.

\begin{figure}[h]
\includegraphics[width=0.3 \textwidth]{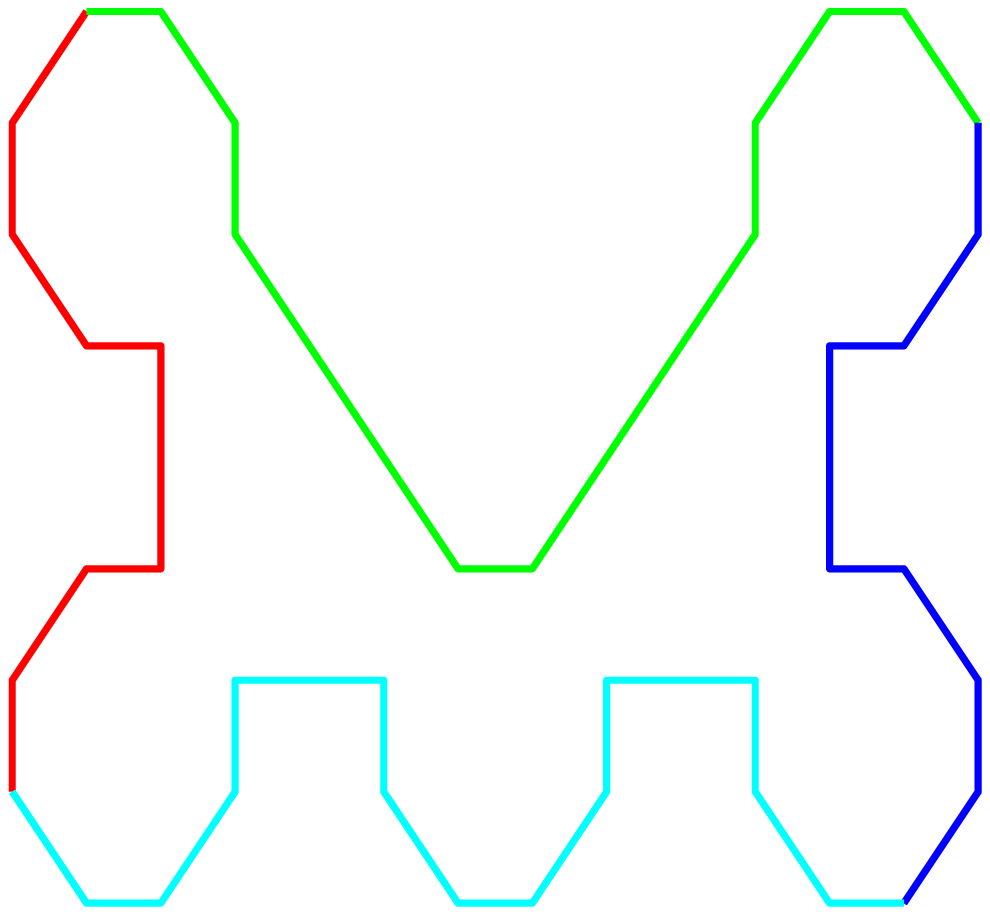}
\includegraphics[width=0.3 \textwidth]{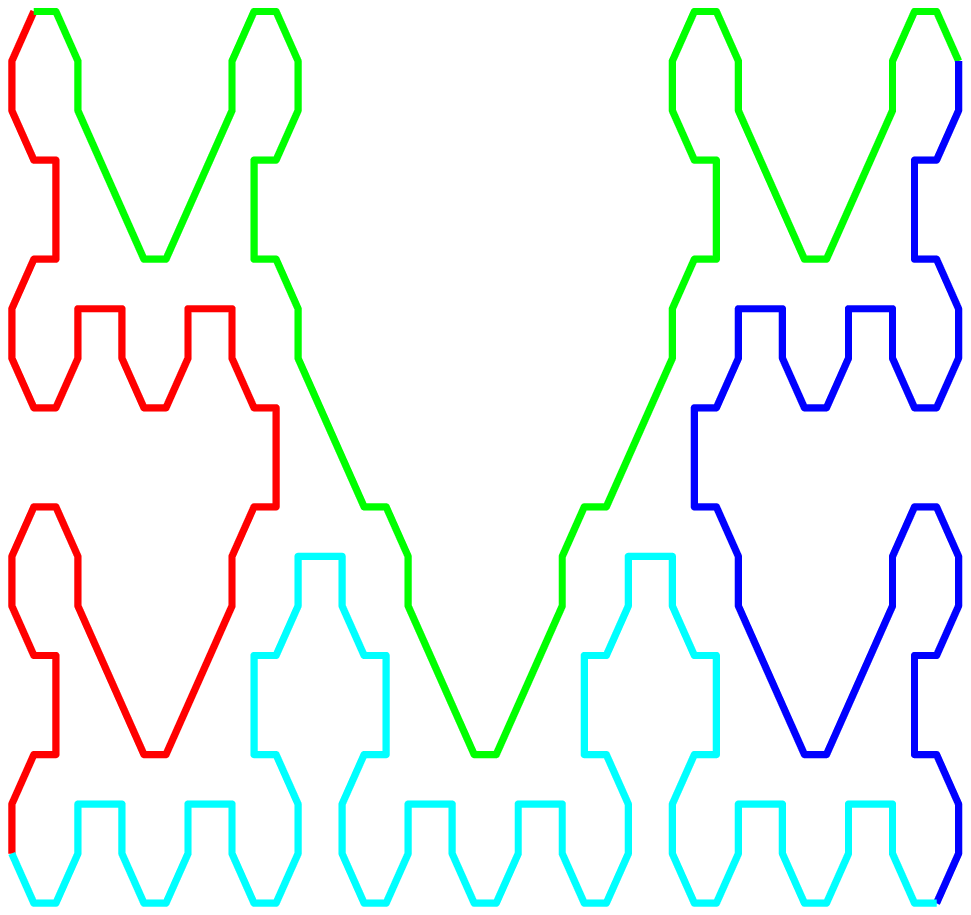}\includegraphics[width=0.3 \textwidth]{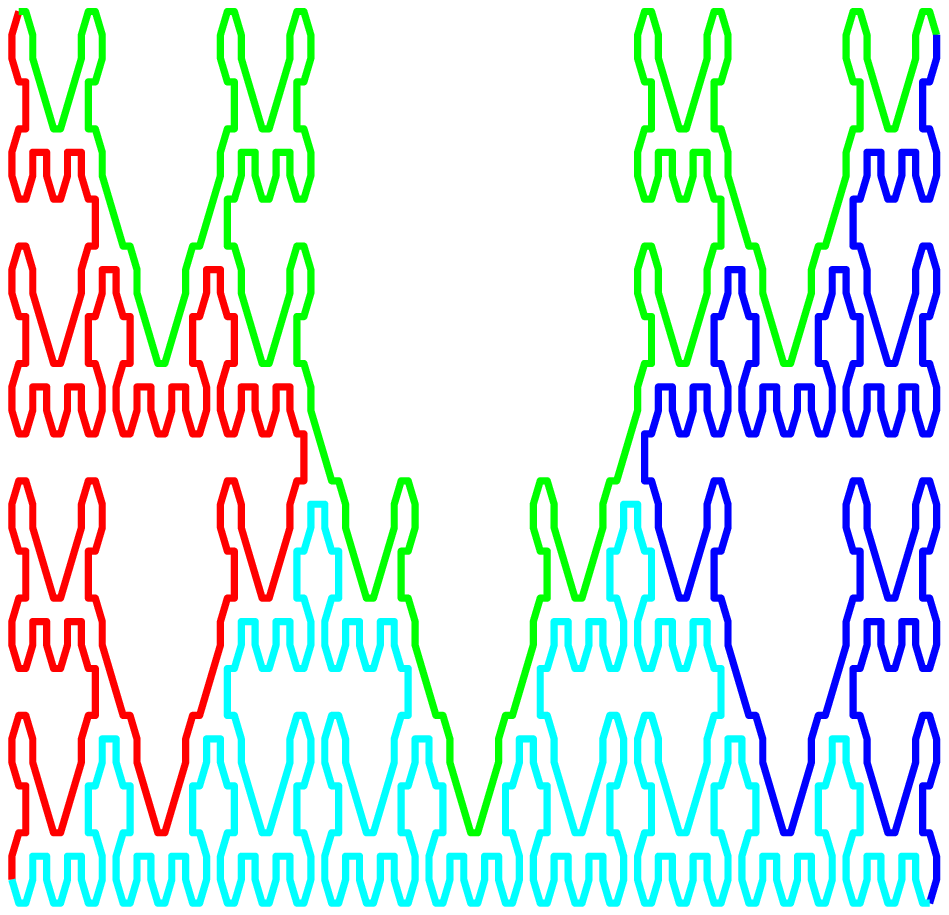}
\caption{The first three approximations to the filling curve of Mcmullen set of Example \ref{Ex2}.}
\label{Ex2_fig3}
\end{figure}

\end{example}

%

\begin{figure}[h]
\includegraphics[width=0.35 \textwidth]{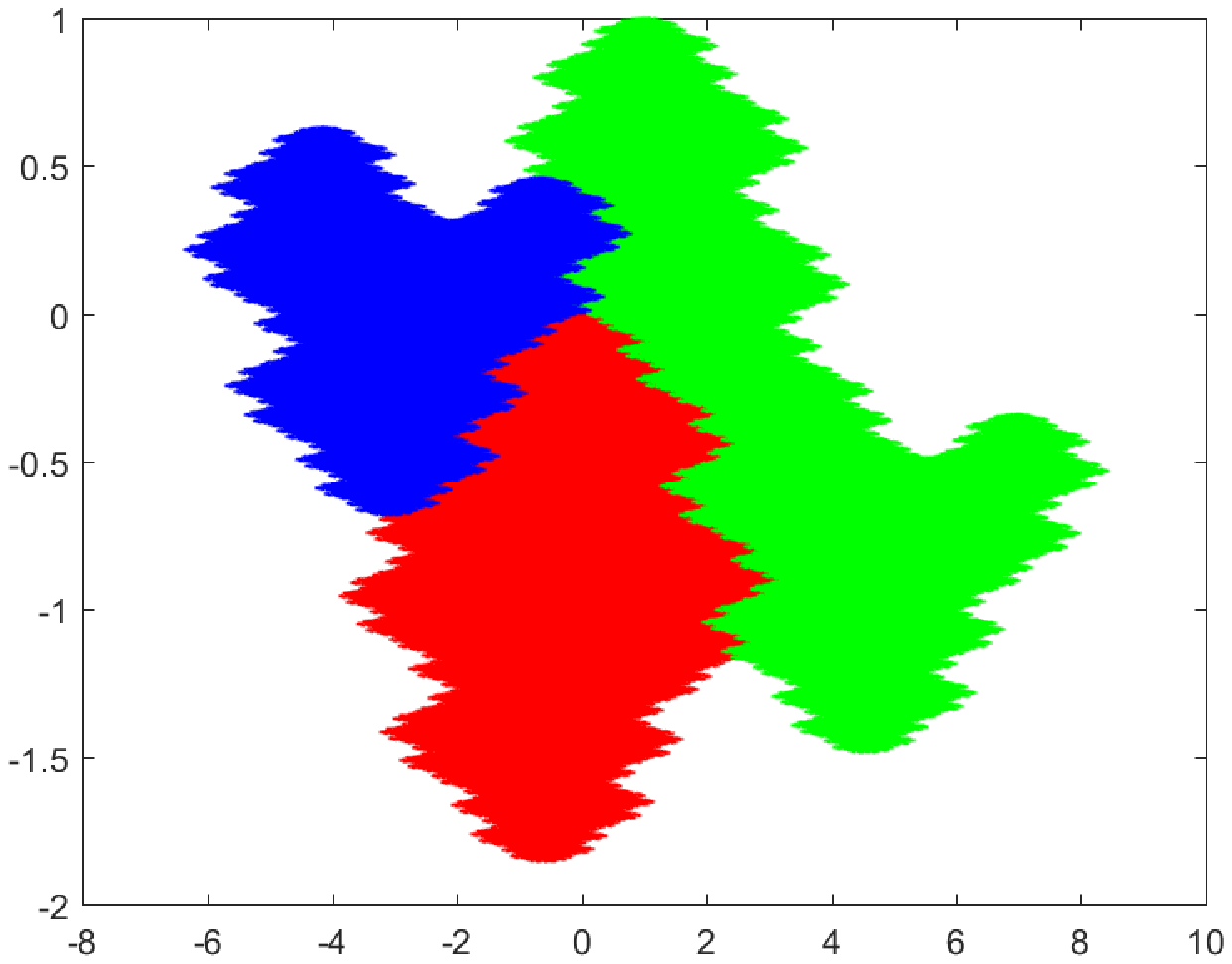}\includegraphics[width=0.35 \textwidth]{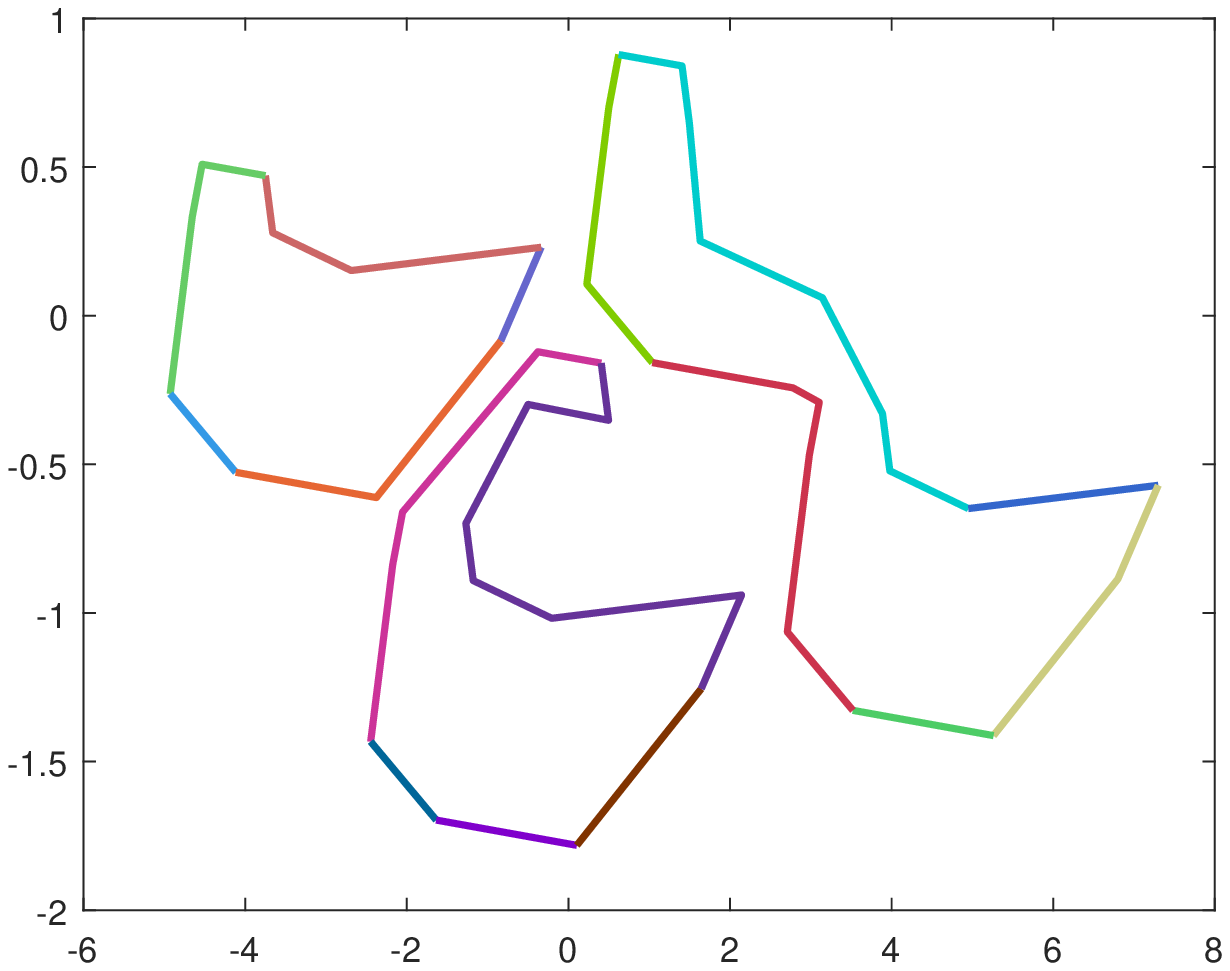}\includegraphics[width=0.35 \textwidth]{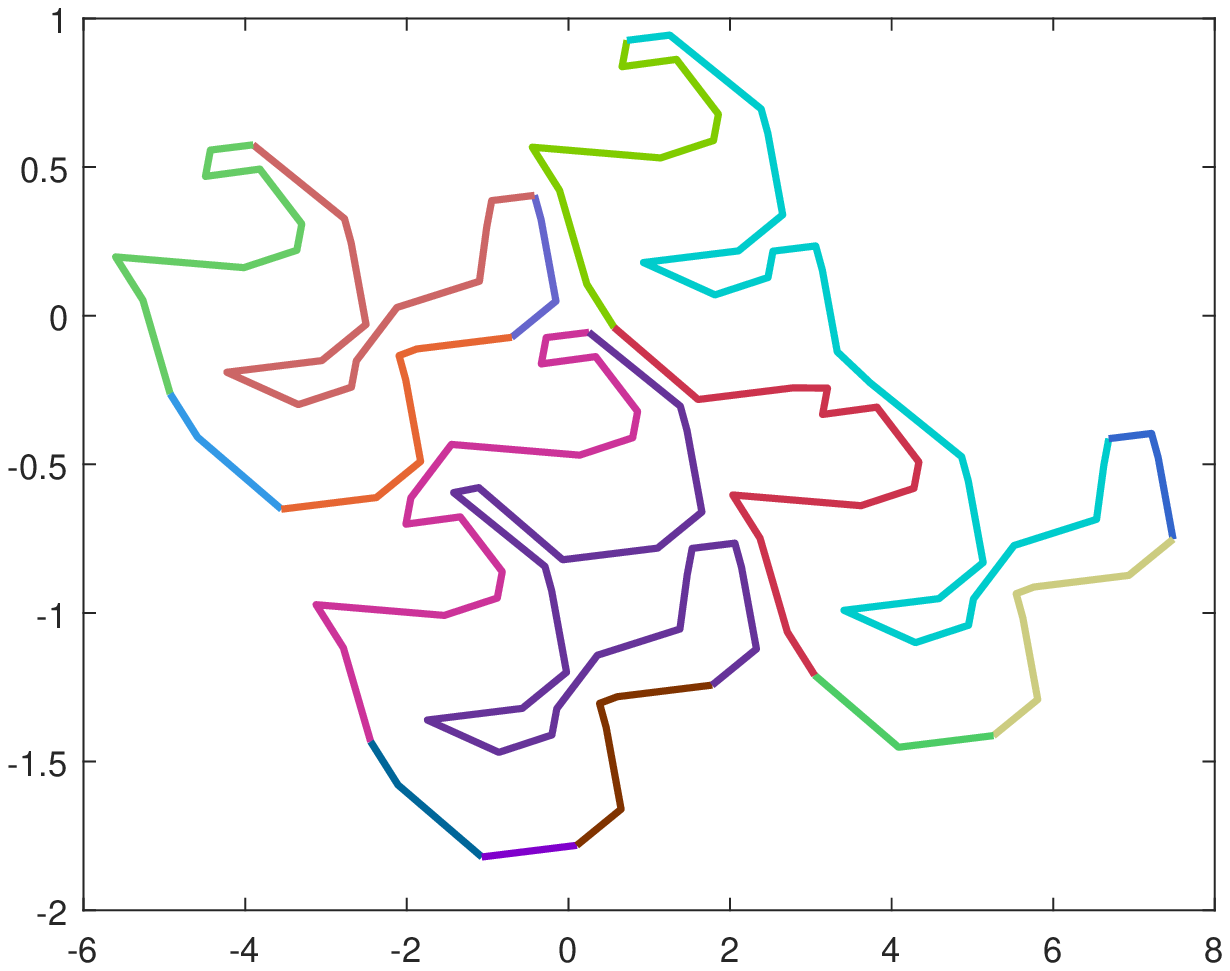}\\
\includegraphics[width=0.4 \textwidth]{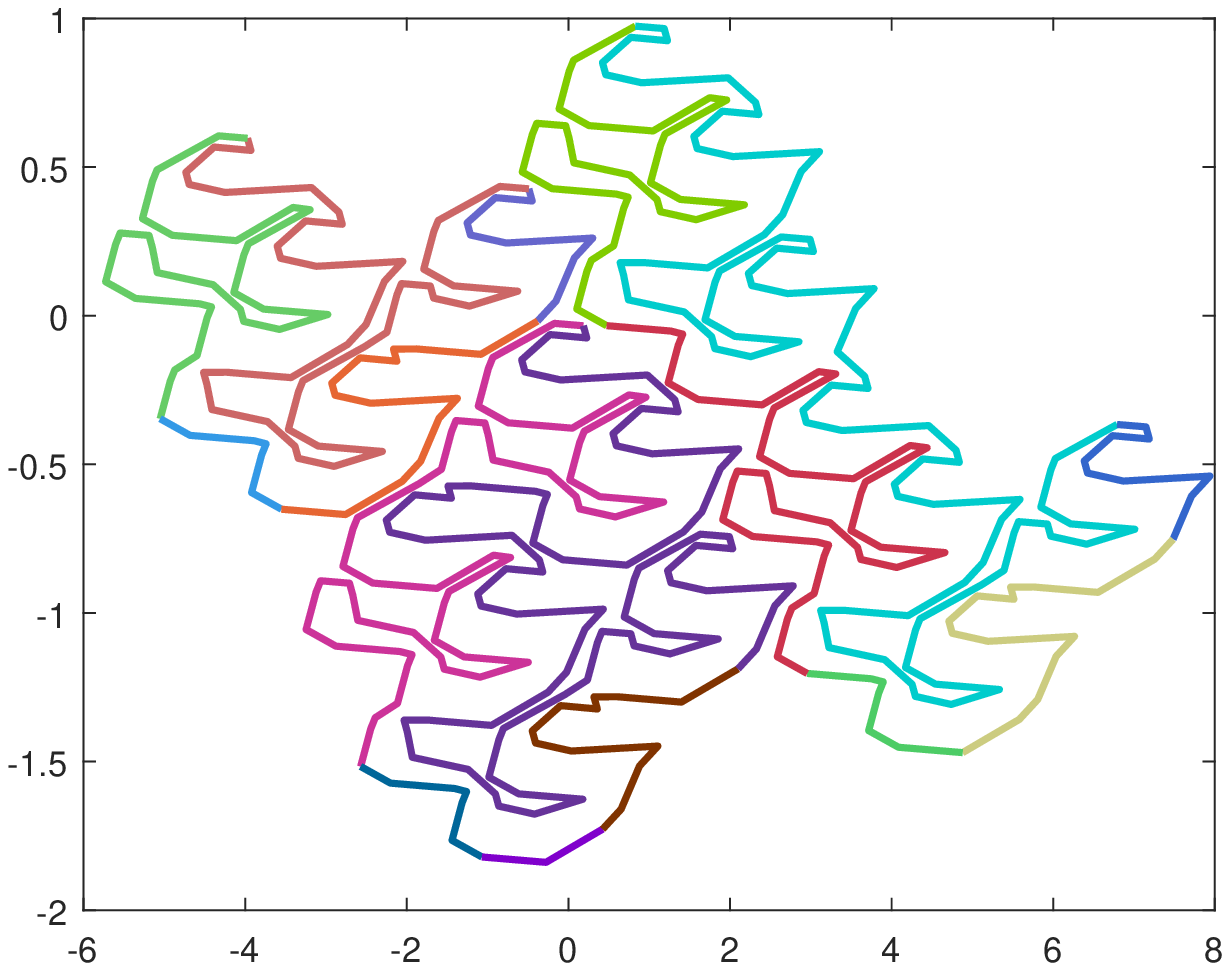}
\includegraphics[width=0.4 \textwidth]{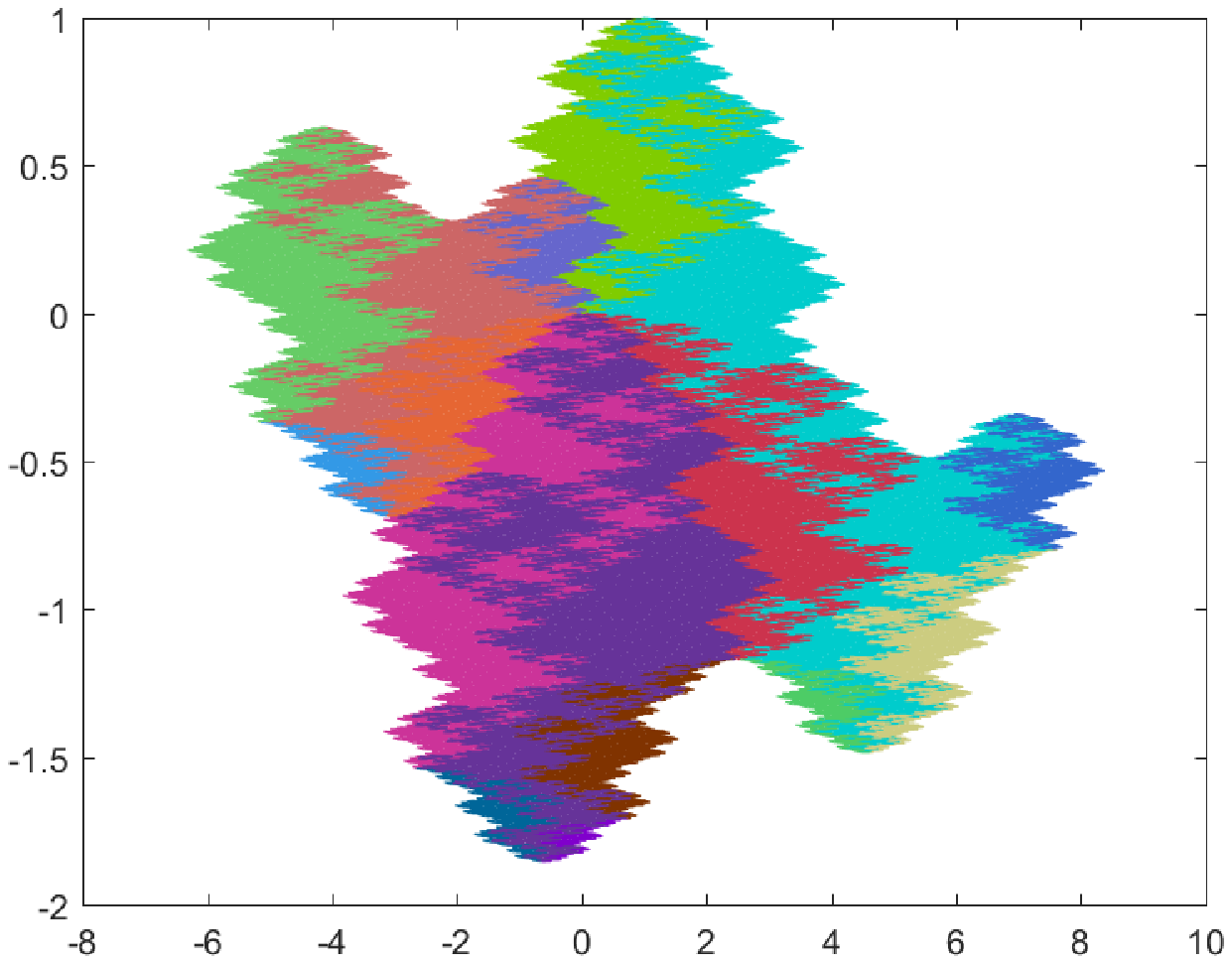}
\caption{The first figure is the Rauzy fractal given by the substition $
\sigma: 1\longrightarrow 12321,2 \longrightarrow 321,3 \longrightarrow 2$. The following three Figures show the first three approximations of the filling curve of this Rauzy fractal.}
\label{Rauzy_fig}
\end{figure}

Inspired by \cite{DaiRaoZhang15}, we shall do some further study on space-filling curve of self-affine fractals, such as space-filling curves of McMullen sets and of Rauzy fractals. To this matter, we try to find a similar method as we did for self-similar sets to construct a linear GIFS from a self-affine set or from a Rauzy fractal (See Figure \ref{Rauzy_fig} where we shows the approximations of the filling curve of a Rauzy fractal). Here we want to emphasis that to do the parametrization of a Rauzy fractal is based on invariant sets of a graph-directed iterated system, that is to say, we construct a linear GIFS from a given GIFS. (\cite{DaiRaoZhang15} focused on constructing a linear GIFS from a given IFS.)
 However, to construct a linear GIFS from the original iterated function system or from graph-directed iterated function system is a complicated work, so we will not do the part of constructions of linear GIFS and only give the set equations of linear GIFS in present paper. Hopefully, it will be shown in our following work. Figure \ref{Rauzy_fig} shows an interesting example of Rauzy fractal.

This paper is organised in the following way. In Section \ref{Preli}, we recall some definitions and notations, and then we give some additional explanations of the Examples~ \ref{Ex3}, \ref{Ex1}, and \ref{Ex2} in Section \ref{Sec_Example}. The last section  is devoted to the verification of our main Theorem \ref{Main1}. 


\section{Preliminaries}\label{Preli}

\subsection{The symbolic space related to a graph $G$} Let $G=(\A,\Gamma)$ be a directed graph.
 A sequence of edges in $G$, denoted by $\bomega=\omega_1+\omega_2+\dots+\omega_n$, is called a \emph{walk}, if the
  terminal state of $\omega_i$ coincides with the initial state of  $\omega_{i+1}$ for $1\leq i \leq n-1$. A walk is called a \emph{trail}, if all the edges appearing in the walk are distinct. A trail is
called a \emph{path} if all the vertices are distinct.
 For $i\in \mathcal{A}$, let
 $$
 \Gamma_i^k,\ \Gamma_i^\ast, \text{ and } \Gamma_{i}^{\infty}
 $$
  be the set of
 all walks of length $k$, the set of all walks of finite length, and the set of all infinite walks,
  starting at state $i$, respectively.
Note that $\Gamma_i^\ast=\bigcup_{k\geq 1}\Gamma_i^k$.

 For $\bomega=(\omega_k)_{k=1}^\infty$, define by $\bomega|_n=\omega_1+\omega_2+\dots+\omega_n$  the prefix of $\bomega$ of length $n$.
Moreover, call
$[\omega_1\dots\omega_n]:=\{{\boldsymbol \gamma}\in \Gamma_i^\infty;~~{\boldsymbol \gamma}|_n=\omega_1+\dots+\omega_n\}$
the \emph{cylinder} associated with a walk $\omega_1+\dots +\omega_n$.

For a walk  ${\boldsymbol \gamma}=\gamma_1+\dots +\gamma_n\in\Gamma_i^n$, set $g_{\bgamma}:=g_{\gamma_1}\circ g_{\gamma_2}\dots\circ g_{\gamma_n}$, then we denote
$$
E_{\boldsymbol \gamma}:=g_{\bgamma}(E_{t(\bgamma)}),
$$
where $t(\bgamma)$ denotes the terminal state of the path $\boldsymbol \gamma$ (which equals$\gamma_n$ here). Iterating \eqref{GIFS1} $k$-times, we obtain
\begin{equation}\label{uni-set3}
E_i=\bigcup_{{\boldsymbol \gamma}\in \Gamma_i^k}E_{\boldsymbol \gamma}.
\end{equation}

We define the projection $\pi: \Gamma_1^{\infty}\times\dots\times\Gamma_N^{\infty} \rightarrow \mathbb{R}^d\times\dots\times\mathbb{R}^d$, where $\pi_i: \Gamma_i^{\infty} \rightarrow \mathbb{R}^d$ is given by
\begin{equation}\label{eq-projection}
\{\pi_i({\boldsymbol \omega})\}:=\bigcap_{n\geq 1}E_{\boldsymbol {\omega|_{n}}}.
\end{equation}
 For $x\in E_i$, we call ${\boldsymbol \omega}$ a \emph{coding} of $x$ if $\pi_i({\boldsymbol \omega})=x$. It is easy to see that $\pi_i(\Gamma_{i}^{\infty})=E_i$.

\subsection{Pseudo-norm and Hausdorff measure in pseudo-norm}

Denote by $B(x,r)$ the open ball with center $x$ and radius $r$. Recall that $A$ is the expanding matrix with $|\det A|=q$, then $V=A(B(0,1))\setminus B(0,1)$ is homeomorphic to an annulus.  For $\delta\in(0,\frac{1}{2})$,  choose a positive $\mathbb{C}^\infty$-function $\phi_{\delta}(x)$  with support in $B(0,\delta)$ such that $\phi_{\delta}(x)=\phi_{\delta}(-x)$ and $\int\phi_{\delta}(x)dx=1$, and then define a pseudo-norm $\|\cdot\|_{\omega}$ in $\R^d$ by
$$\|x\|_{\omega}=\sum_{n\in \Z}q^{-n/d}h(A^{n}x),$$
where $h(x)=\chi_V*\phi_{\delta}(x)=\int_{\R^d}\chi_V(x-y)\cdot\phi_{\delta}(y)dy$.

We list some basic properties of $\|\cdot\|_{\omega}$.

\begin{proposition}(See \cite[Proposition 2.1]{HeLau08}) \label{Pro-Weak}
The function $\|\cdot\|_{\omega}$ satisfies the following conditions.

(i) $\|x\|_{\omega}\geq 0$; $\|x\|_{\omega}=0$ if and only if $x=0$.

(ii) $\|x\|_{\omega}=\|-x\|_{\omega}$;

(iii) $\|Ax\|_{\omega}=q^{1/d}\|x\|_{\omega}\geq\|x\|_{\omega}$ for all $x\in \R^d$.

(iv) There exists  a constant $\beta>0$ such that $\|x+y\|_{\omega}\leq \beta \max \{\|x\|_{\omega}, \|y\|_{\omega}\}$ for any $x,y \in \R^d$.

\end{proposition}

The pseudo-norm $\|\cdot\|_{\omega}$ is comparable with the Euclidean norm $\|x\|$ .

\begin{proposition}\label{weaknorm}
(See \cite[Proposition 2.4]{HeLau08})  Let  $\lambda_{\max}$ and  $\lambda_{\min}$  be  the maximal modulus and minimal modulus of the eigenvalues of $A$, respectively. For any $0<\varepsilon<\lambda_{\min}-1$, there exists $C>0$(depends on $\varepsilon$) such that
$$ C^{-1}\|x\|^{\ln q/d \ln(\lambda_{\max}+\varepsilon)}\leq \|x\|_{\omega}\leq C\|x\|^{\ln q/d \ln(\lambda_{\min}-\varepsilon)},\quad\text{ if } \|x\|>1;$$
$$C^{-1}\|x\|^{\ln q/d \ln(\lambda_{\min}-\varepsilon)}\leq \|x\|_{\omega}\leq C\|x\|^{\ln q/d \ln(\lambda_{\max}+\varepsilon)},\quad\text{ if } \|x\|\leq 1.$$
\end{proposition}

The Hausdorff measure  with respect to the pseudo-norm $\|\cdot\|_{\omega}$  was given by He and Lau \cite{HeLau08} as follows.
For $E\subset\R^d$, set $\text{diam}_{\omega}E=\text{sup}\{\|x-y\|_{\omega};~~x,y\in E\}$ to be the $\omega$-diameter of $E$. For $s\geq 0, \delta>0$, set
$$\mathcal{H}_{\omega,\delta}^s(E)=\inf\big\{\sum_{i=1}^\infty (\text{diam}_{\omega}E_i)^s;~~E\subset \bigcup_i E_i,~~ \text{diam}_{\omega}E_i\leq\delta\big\},$$
$$\mathcal{H}_{\omega}^s(E)=\lim_{\delta\rightarrow 0}\mathcal{H}_{\omega,\delta}^s(E).$$
$\H_{\omega}^s$ has the translation invariance property and the scaling property \cite{HeLau08}, that is,
$$\H_\omega^s(E+x)=\H_\omega^s(E) \text{ and } \H_{\omega}^s(A^{-1}E)=q^{-s/d}\H_{\omega}^{s}(E).$$
Thus the Hausdorff dimension with respect to $\|\cdot\|_{\omega}$ can be defined by
$$\dim_{\omega} E=\inf\{s;~~\H_{\omega}^s(E)=0\}=\sup\{s;~~\H_{\omega}^s(E)=\infty\}.$$

\subsection{Linear GIFS}
Let $({\mathcal{A}},\Gamma,\mathcal{G})$ be a GIFS with vertex set $\mathcal{A}$, edge set $\Gamma$ and mapping set $\mathcal{G}$.
Let $\Gamma_i=\Gamma_i^1$ be the set of outgoing edges from the state $i$.
We call $({\mathcal{A}},\Gamma, \mathcal{G}, \prec)$ an \emph{ordered GIFS},
if  $\prec$ is a partial order on $\Gamma$ such that
\begin{enumerate}
  \item[($i$)] $\prec$ is a linear order when restricted on $\Gamma_j$ for every $j\in {\mathcal A}$;
  \item[($ii$)]  elements in $\Gamma_i \text{ and } \Gamma_j$ are not comparable if  $i \neq j$.
\end{enumerate}
(See \cite{RaoZhangS16} for detail.)

 The  order $\prec$  induces a lexicographical order on each
 $\Gamma_i^k$.
  Observe that $(\Gamma_i^k, \prec)$ is a linear order;
  two paths  ${\boldsymbol \gamma},{\boldsymbol \omega}\in \Gamma_{i}^k$ are said to be \emph{adjacent}
  if there is no walk between them  with respect to the order $\prec$.

\begin{definition}(see \cite{RaoZhangS16})\label{Linear} {\rm An ordered GIFS $({\mathcal{A}},\Gamma,\mathcal{G},\prec)$   with invariant sets $\{E_i\}_{i=1}^N$
 is called a \emph{linear} GIFS,  if for all $i\in {\mathcal A}$ and $k\geq 1$,
 $$
 E_{\boldsymbol \gamma}\cap E_{\boldsymbol \omega}\neq \emptyset
 $$
  for adjacent walks ${\boldsymbol \gamma},{\boldsymbol \omega}$  in $\Gamma_{i}^k$ .}
\end{definition}
For $i\in \mathcal{A}$, a walk $\bomega\in\Gamma_i^{\infty}$ is called the \emph{lowest} walk, if $\bomega|_n$ is the lowest walk in $\Gamma_i^n$ for all $n$; in this case, we call $a=\pi_i(\bomega)$ the \emph{head} of $E_i$.
Similarly, we define the highest walk $\bomega'$ of $\Gamma_i^\infty$, and we call $b=\pi_i(\bomega')$ the
  the \emph{tail} of $E_i$.

\begin{definition}(see \cite[Definition 4.1]{RaoZhangS16} )\label{Chain}{\rm
An ordered GIFS is said to satisfy the \emph{chain condition}, if for any $i\in {\mathcal A}$, and any two adjacent edges
$\omega, \gamma\in \Gamma_i$ with $\omega \prec \gamma$,
$$g_\omega(\text{tail of } E_{t(\omega)})=g_\gamma(\text{ head of }E_{t(\gamma)}).$$
}
\end{definition}
\begin{lemma}\label{Text_linear}
An ordered GIFS is a linear GIFS if and only if it satisfies the chain condition.
\end{lemma}
{\bf Remark.} Definition \ref{Linear}, Definition \ref{Chain} and Lemma \ref{Text_linear} still make sense when  $({\mathcal{A}},\Gamma,\mathcal{G})$ is a single matrix GIFS.

\section{Examples}\label{Sec_Example}
In this section, we will give the set equation form of substitution rules which we showed in the previous examples and explain simply that these induced GIFS are linear.
\begin{example}\label{Ex31}
{\bf A unit square.} Through the edge-to-trail substitution rule, we get the following GIFS in the set equation form.
\begin{flalign*}
  E_1\circeq\vphantom{a}&S_1(E_1)\cup S_1(E_2)\cup S_3(E_1)\cup S_4(E_3)\cup S_4(E_1)\cup S_5(E_1),
\\
  E_2\circeq\vphantom{a}& S_5(E_2)\cup S_5(E_3)\cup S_4(E_2)\cup S_6(E_1)\cup S_6(E_2),\\
  E_3\circeq\vphantom{a}& S_6(E_3)\cup S_3(E_2)\cup S_3(E_3)\cup S_2(E_2)\cup S_2(E_3)\cup S_2(E_1)\cup S_1(E_3).
\end{flalign*}
Here we use $``\circeq"$ to emphasize that it is an ordered GIFS. 
By the result of \cite{DaiRaoZhang15}, the above induced GIFS is a linear GIFS, or we can check it by Lemma \ref{Chain} directly.

\end{example}

\begin{example}\label{Ex11}
\begin{figure}[h]
\includegraphics[width=0.4 \textwidth]{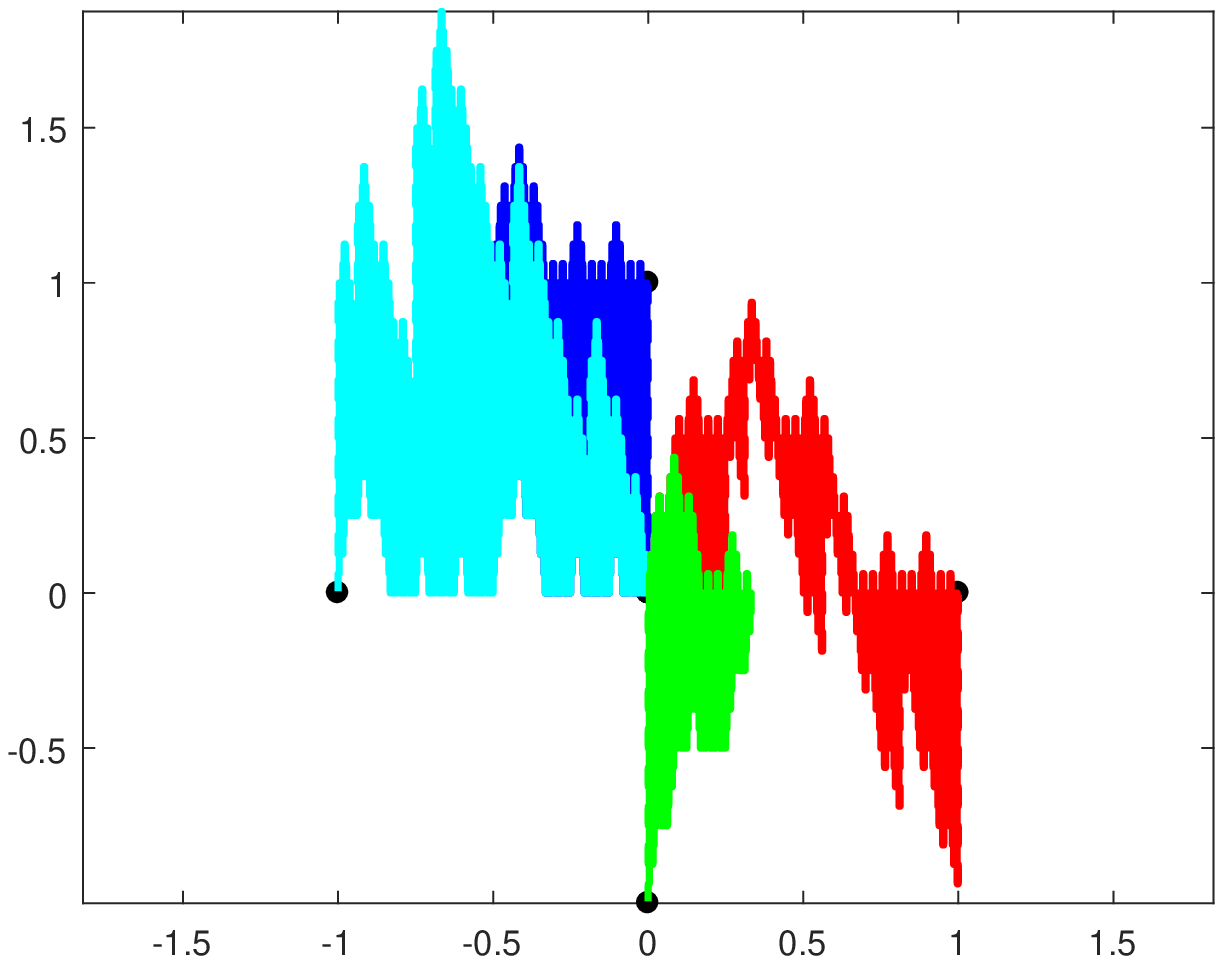}\includegraphics[width=0.4 \textwidth]{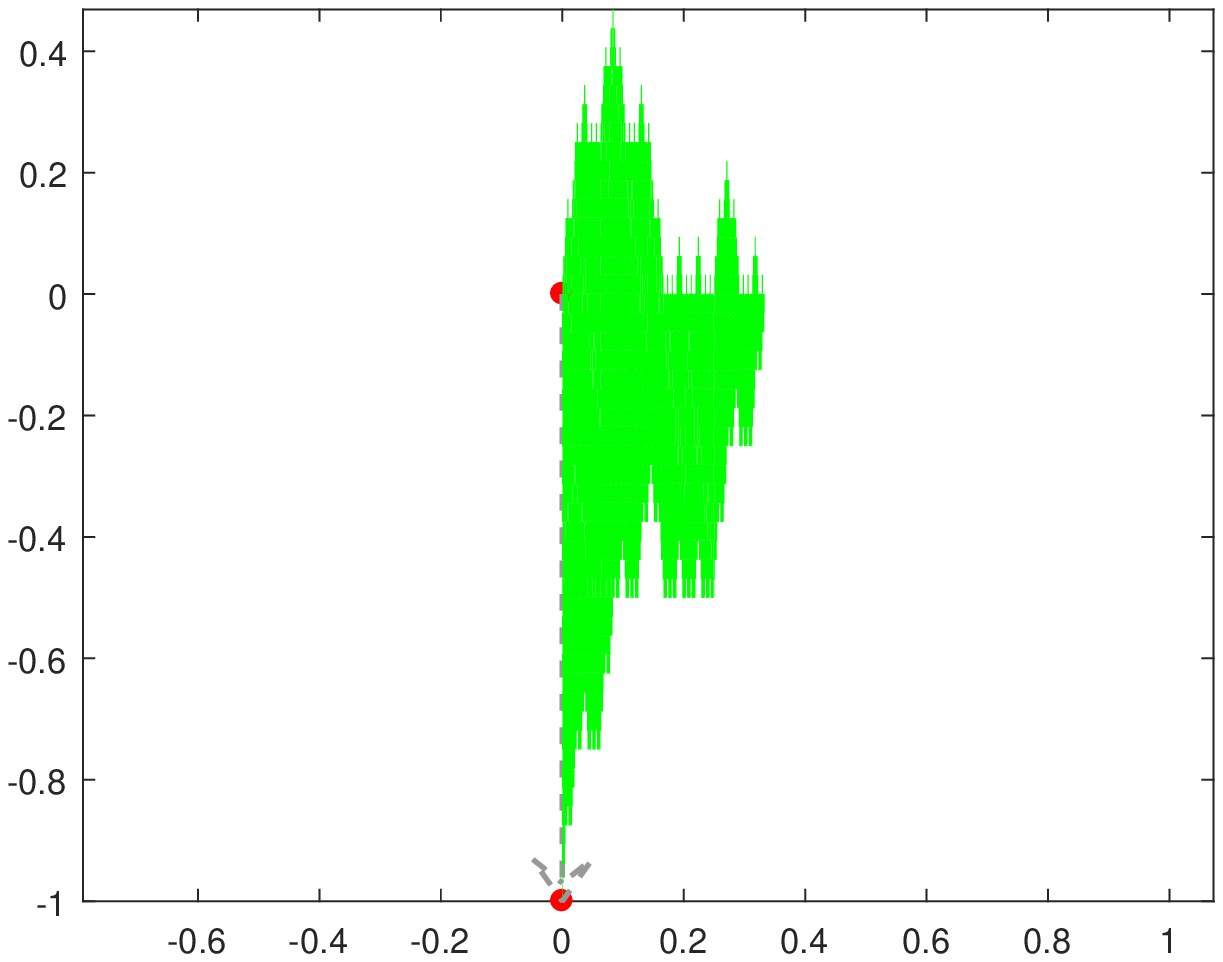}
\caption{The left figure is $E=\cup_{i=1}^4E_i$, and the right one is $E_4$. Example \ref{Ex1}.}
\end{figure}

 {\bf Dekking's plane filling curve}   \cite{Dekking82b}.
Here we continue Example \ref{Ex1}. We translate the substitution of Dekking's example to the following set equation
{\small
\begin{flalign*}
AE_1\circeq\vphantom{a}& E_1\cup(E_2+\e_1)\cup(E_1+\e_1+\e_2)\cup(E_4+2\e_1+\e_2)
\cup(E_1+2\e_1)\cup\\ &(E_4+3\e_1)
\cup(E_1+3\e_1-\e_2)\cup(E_2+4\e_1-\e_2),\\
AE_2\circeq\vphantom{a}& E_3\cup(E_2-\e_1)\cup(E_3-\e_1+\e_2)\cup(E_2-2\e_1+\e_2)
\cup(E_1-2\e_1+2\e_2)\cup\\
&(E_4-\e_1+2\e_2)\cup(E_1-\e_1+\e_2) \cup(E_2+\e_2),\\
AE_3\circeq\vphantom{a}& E_3\cup(E_2-\e_1)\cup(E_3-\e_1+\e_2)\cup(E_2-2\e_1+\e_2)\cup(E_3-2\e_1+2\e_2)\cup
\\ &(E_4-3\e_1+2\e_2)\cup(E_1-3\e_1+\e_2)\cup(E_4-2\e_1+\e_2)\cup(E_3-2\e_1)\cup
\\ &(E_2-3\e_1)\cup(E_3-3\e_1+\e_2)\cup(E_4-4\e_1+\e_2),
\\
AE_4\circeq\vphantom{a}& E_1\cup(E_4+\e_1)\cup(E_3+\e_1-\e_2)\cup(E_4-\e_2).
\end{flalign*}
}
To check the chain condition, we need to calculate the heads and tails of $E_1, E_2, E_3, E_4$. We denote the head of a set $K$ by $h(E)$ and tail of a set $t(E)$. Then we have
 \begin{align*}
 h(E_1)=0,~t(E_1)=\e_1,& ~h(E_2)=0, ~t(E_2)=\e_2,\\
h(E_3)=0,~ t(E_3)=-\e_1,& ~h(E_4)=0,~t(E_4)=-\e_2.
\end{align*}
Thus it is easy to check that it satisfies the chain condition. Figure \ref{Ex1_Fig2} shows the proceeding of filling curve of $E_1$.
\end{example}
\begin{example}\label{Ex2plus}

{\bf A McMullen set}.
The McMullen set $T$ (see Figure \ref{Ex2_fig1}, left) in Example \ref{Ex2} is given by $T=\bigcup_{i=1}^5 S_i(T)$ with $$\Big\{S_i(x)=A^{-1}(x+d_{i})\Big\}_{i=1}^5,\quad d_1=0, d_2=\e_1, d_3=2\e_1, d_4=\e_2, d_5=\e_1+\e_2.$$
Through the substitution rule in Example \ref{Ex2}, we obtain the following set equation form of a GIFS.
\begin{flalign*}
  E_1\circeq\vphantom{a}&S_1(E_1)\cup S_1(E_2)\cup S_4(E_4)\cup S_4(E_1),
\\
  E_2\circeq\vphantom{a}& S_4(E_2)\cup S_4(E_3)\cup S_2(E_2)\cup S_5(E_1)\cup S_5(E_2),\\
  E_3\circeq\vphantom{a}& S_5(E_3)\cup S_5(E_4)\cup S_3(E_2)\cup S_3(E_3),\\
 E_4\circeq\vphantom{a}&S_3(E_4)\cup S_3(E_1)\cup S_2(E_3)\cup S_2(E_4)\cup S_2(E_1)\cup S_1(E_3)\cup S_1(E_4).
\end{flalign*}
In the same way as Example \ref{Ex11}, we check that the above GIFS satisfied the chain condition. Then it is a linear GIFS and clearly the open set condition is satisfied. Actually the union of the invariant sets $\bigcup_{i=1}^4 E_i$ is the McMullen set $T$.
\end{example}

\section{Proof of Theorem \ref{Main1}}\label{Sec_Main1}
In this section, 
 we prove Theorem \ref{Main1} by constructing an auxiliary GIFS (which we call measuring-recording GIFS), which is very similar to
 the proof in \cite{RaoZhangS16}.
However, the theorem related to the open set condition of Mauldin and Williams \cite{MauldinWilliams88} does  not hold when $A$ is not a similitude.
So we need to use the result of Luo and Yang \cite{LuoYang10} to modify the proof.



\subsection{Markov measure induced by GIFS}
Let $(\mathcal{A},\Gamma,\mathcal{G})$ be single-matrix GIFS with expanding matrix $A$ and  $\{E_i\}_{i=1}^N$ be the invariant sets. Denote $q=|\text{det} (A)|$.
And $M=(m_{ij})_{1\leq i,~j\leq N}$ is the associated matrix of the directed graph $(\mathcal{A},\Gamma)$. Due to the following lemma from \cite{LuoYang10}, we can construct the Markov measure.

  \begin{lemma}\label{YL} (\cite[Theorem 1.2]{LuoYang10})
  For  a single matrix GIFS $(\mathcal{A},\Gamma,\mathcal{G})$,  let $\lambda$ be the maximal eigenvalue of $M$.  If $M$ is primitive and the OSC holds, then for any $1\leq i \leq N$,

  (i) $\alpha=\dim_\omega E_i=d\frac{\log\lambda}{\log q}$;

  (ii) $0<\H_\omega^{\alpha}(E_i)<\infty$.

  (iii)  The right hand side of \eqref{GIFS1} is a disjoint union in sense of the measure of $\H_\omega^{\alpha}$.
\end{lemma}
{\bf Remark $(1)$.}  By item $(iii)$ of the above lemma, we immediately have
 $${\mathcal H}_\omega^{\alpha}(E_{\boldsymbol \omega}\cap E_{\boldsymbol \gamma})=0$$
  for any incomparable ${\boldsymbol \omega},{\boldsymbol \gamma} \in \Gamma_i^\ast$. (Two walks are said to be \emph{comparable} if one of them is a prefix of the other.)

{\bf Remark $(2)$.} Since $E_i=\bigcup_{j=1}^N\bigcup_{e\in \Gamma_{ij}}g_e(E_j)$, using Remark (1),  we get
  $$\H_\omega^{\alpha}(E_i)=\sum_{j=1}^{N}\sum_{e\in \Gamma_{ij}}\H_\omega^{\alpha}(g_e(E_j))=\lambda^{-1}\sum_{j=1}^{N}\sharp\Gamma_{ij}\H^{\alpha}_\omega(E_j).$$
This shows that $(\mathcal{H}_{\omega}^{\alpha}(E_1),\dots,\mathcal{H}_{\omega}^{\alpha}(E_N))$ is an eigenvector with respect to $\lambda$ of $M$.

In the rest of the section, we will always assume that ${\cal G}$ satisfies the conditions of Lemma \ref{YL}. Then $0<\H_\omega^{\alpha}(E_i)<\infty$ for all $1\leq i \leq N$.
Now, we define Markov measures on the symbolic spaces $\Gamma_i^\infty$, $i\in\mathcal{A}$.
For arbitrary edge $e\in \Gamma$ such that $e\in \Gamma_{ij}$, set
\begin{equation} \label{our-weight}
p_e=\frac{\H_{\omega}^{\alpha}(E_j)}{\H_{\omega}^{\alpha}(E_i)}\lambda^{-1}.
\end{equation}
Using Remark $(2)$ of Lemma \ref{YL},
it is easy to verify that $(p_e)_ {e\in \Gamma}$ satisfies
\begin{equation}\label{weights}
\sum_{j\in {\mathcal A}}\sum_{e\in \Gamma_{ij}} p_e=1, \text{ for all } i\in \mathcal{A}.
\end{equation}
We  call $(p_e)_ {e\in \Gamma}$ a \emph{probability weight vector}.
 Let $\mathbb{P}_i$ be a Borel measure on $\Gamma_i^\infty$ satisfying the relations
\begin{equation}\label{measure}
\mathbb{P}_i([\omega_1\dots\omega_n])=\H_{\omega}^{\alpha}(E_i) p_{\omega_1}\dots p_{\omega_n}
\end{equation}
for all cylinder $[\omega_1\dots\omega_n]$. The existence of such measures is guaranteed by \eqref{weights}. We call $\{{\mathbb P_i}\}_{i=1}^N$ the \emph{Markov measures} induced by the GIFS ${\mathcal G}$.

Denote the restriction of ${\mathcal H}^{\alpha}_{\omega}$ on $E_i$ by $\mu_i={\mathcal H}^{\alpha}_{\omega}|_{E_i} , \text{ for }i=1,\dots, N.$ The following Lemma gives the relation between the Markov measure and the restricted Hausdorff measure.

 \begin{lemma}(see \cite{MauldinWilliams88, LuoYang10})\label{lem-Markov}
Suppose the single-matrix graph IFS $(\A, \Gamma, \G )$ satisfies the OSC and
  the associated matrix $M$ is primitive.
 Let $\pi_i: \Gamma_i^\infty\to E_i$ be the projections defined by \eqref{eq-projection}. Then
$$
\mu_i=\mathbb{P}_i\circ \pi_i^{-1}.
$$
 \end{lemma}



\subsection{The construction of measure-recording GIFS}

Let $(\mathcal{A},\Gamma,\mathcal{G},\prec)$ be a linear GIFS such that the open set condition is fulfilled and the associated matrix is primitive, then
$0<{\cal H}_\omega^{\alpha}(E_i)<\infty$ for all $i$, where $\alpha=d\frac{\log\lambda}{\log q}$ by Lemma \ref{YL}.

For $i\in \mathcal{A}$, we list the edges in $\Gamma_i$ in the ascendent order with respect to $\prec$, i.e.,
$$
\gamma_1\prec \gamma_2\prec\dots\prec \gamma_{\ell_i}.
$$
 Recall that $t(\gamma)$ denotes the terminate vertex of an edge $\gamma$. Then by \eqref{GIFS1}, we can rewrite $E_i$ as
$$
E_i\circeq g_{\gamma_1}(E_{t(\gamma_1)})\cup \cdots \cup g_{\gamma_{\ell_i}}(E_{t(\gamma_{\ell_i})}).
$$
Here we use `$\circeq$' to emphasize the order of the union of the right side.

Denote by
$F_i=[0, {\cal H}_\omega^{\alpha}(E_i)]$ an interval on $\mathbb{R}$, then by equation \eqref{weights}, we have
\begin{equation}\label{F_i}
F_i=\Big[0,~~\H_{\omega}^{\alpha}(g_{\gamma_1}(E_{t(\gamma_1)}))\Big]\cup\dots\cup\Big[\sum_{j=1}^{{\ell}_i-1}\H_{\omega}^{\alpha}(g_{\gamma_j}(E_{t(\gamma_j)})), ~~\sum_{j=1}^{{\ell}_i}\H_{\omega}^{\alpha}(g_{\gamma_j}(E_{t(\gamma_j)})) \Big].
\end{equation}
We define the mappings,
$$
f_{\gamma_k}(x)=q^{-\alpha/d}x+b_k: ~~\mathbb{R}\longrightarrow\mathbb{R}, \quad 1\leq k \leq \ell_i,
$$
where $b_k=\sum_{j=1}^{k-1}\H_{\omega}^{\alpha}(E_{t(\gamma_j)})q^{-\alpha/d}$. Then $F_i$ satisfies the following equation by \eqref{F_i}
\begin{equation}\label{partition-F}
F_i\circeq f_{\gamma_1}(F_{t(\gamma_1)})\cup \cdots \cup f_{\gamma_{\ell_i}}(F_{t(\gamma_{\ell_i})}).
\end{equation}
Repeating these procedures for all $i\in {\mathcal A}$, equation \eqref{partition-F} gives us an ordered GIFS on $\mathbb{R}$. Set $\mathcal{F}=\{f_{\gamma}: \mathbb{R}\longrightarrow\mathbb{R};~~ \gamma\in\Gamma\},$
and denote this GIFS by
$$
(\mathcal{A}, \Gamma, {\cal F}, \prec),
$$
and call it the \emph{measure-recording GIFS} of $({\mathcal A}, \Gamma, {\mathcal G}, \prec)$. And the invariant sets of the measure-recording  GIFS are $\{F_i\}_{i=1}^N$. (See \cite{RaoZhangS16}.)

Obviously, the measure-recording GIFS has the same graph and the same order as the original GIFS; also keeps the Hausdorff measure information of the original GIFS.  And it is easy to check ${\mathcal F}$ satisfies the open set condition. In fact, the open intervals $\{U_i=(0,\H_{\omega}^{\alpha}(E_i))\}_{i=1}^N$ are the according open sets.

For an edge $e\in \Gamma$, the contraction ratio of $f_e$ is $q^{-\alpha/d}=\lambda^{-1}$, then it is easy to check $(\mathcal{L}(F_1),\dots,\mathcal{L}(F_N))$ is an eigenvector of $M$ with respect the eigenvalue $\lambda^{-1}$. Thus the Markov measure induced by the measure-recording GIFS coincides with $\{\mathbb{P}_i\}_{i=1}^N$ induced by the original GIFS.

Let
$$
\pi_i: \Gamma_i^\infty \rightarrow E_i \text{ and }
\rho_i: \Gamma_i^\infty \rightarrow F_i,\quad i=1,\dots, N,
$$
 be  projections w.r.t. the GIFS $({\cal G})$ and $({\cal F})$, respectively, (see \eqref{eq-projection}). Define
\begin{equation}\label{psi-map}
\psi_i:=\pi_i\circ \rho_i^{-1}.
\end{equation}
In \cite{RaoZhangS16}, it is shown that $\psi_i$ is a well-defined mapping from $F_i$ to $E_i$ since we consider a linear GIFS.



%

Now, we prove Theorem \ref{Main1} by showing that the mapping $\psi_i$ is an optimal parametrization of $E_i$.
\medskip

\noindent {\bf Proof of Theorem \ref{Main1}.}
We use the same notations as before,
let ${\cal F}$ be the measure-recording GIFS of ${\cal G}$. Through the discussion before, we denote the common Markov measure
induced by ${\mathcal G}$ and ${\mathcal F}$ by ${\mathbb P}_i$. $\psi_i=\pi_i\circ\rho_i^{-1}$ is the well-defined mapping from $F_i$ to $E_i$. Let
 $\nu_i={\mathcal L}|_{F_i}$ be the restriction of the Lebesgue measure on $F_i$ and $\mu_i={\mathcal H}_{\omega}^{\alpha}|_{E_i}$ be the restriction of the weak Hausdorff measure on $E_i$, then  $\nu_i=\mathbb{P}_i\circ\rho_i^{-1}$, $\mu_i=\mathbb{P}_i\circ\pi_i^{-1}$ by Lemma \ref{lem-Markov}.

The fact that $\psi_i$ is almost one to one and measure preserving follows by the same arguments as in the self-similar case and we refer to the proof of \cite[Theorem 1.1]{RaoZhangS16}.

\begin{figure}[h]
  \centering
  \includegraphics[width=0.55 \textwidth]{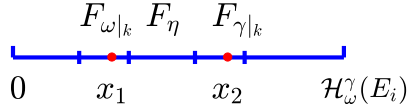}
 \vspace{-1.0 cm}  \caption{ ${\bomega}$ and ${\bgamma}$ are not adjacent.} \label{unadjacent}
\end{figure}

We have to prove the $1/\alpha$-H\"older continuity of $\psi_i$.  From the previous construction, we know that $F_i=[0,\H_{\omega}^{\alpha}(E_i)]$. Now we choose two different points $x_1, x_2$ from $F_i$ which are determined by $\bomega=(\omega_i)_{i=1}^\infty$ and $\bgamma=(\gamma_i)_{i=1}^\infty$, respectively, that is,
$x_1=\rho_i(\bomega), x_2=\rho_i(\bgamma)$. Let $k$ be the smallest integer such that $x_1, x_2$ belongs to two different cylinders. Set $\bomega|_k=\omega_1+\dots+\omega_k$, we know that $\bgamma|_k$ is only different from $\bomega|_k$ at last edge, i.e., $\bgamma|_k=\omega_1+\dots+\omega_{k-1}+\gamma_k$.
%
%
We consider two cases according to ${\bomega}|_k$ and ${\bgamma}|_k$ are adjacent or not.
First, we consider that ${\omega}_k$ and ${\gamma}_k$ are not adjacent. (See Figure \ref{unadjacent}.)
 Then there is a cylinder
$\bdeta=\omega_1+\dots +\omega_{k-1}+\eta_k$ between $\bomega|_k$ and $\bgamma|_k$, so
$$
\|x_1-x_2\|\geq \text{diam}~F_{\boldsymbol \eta}\geq
 h  \cdot (q^{-\alpha/d})^k,
$$
where $h=\min\{{\cal H}_\omega^{\alpha}(E_i);~i=1,\dots, N\}$.
Since $x_1$ and $x_2$ belong to $\rho([\omega_1\omega_2\dots\omega_{k-1}])$ and denote $\bomega_*=\omega_1+\dots +\omega_{k-1}$, the images of $x_1$ and $x_2$ under $\pi_i\circ\rho_i^{-1}$, which denote by $y_1$ and $y_2$, respectively,
belong to $\pi_i([{\bomega}_*])=E_{{\bomega}_*}$.
Then we have
\begin{equation}\label{Holder-1}
\begin{split}
\|y_1-y_2\|_{\omega}
&\leq \text{diam}_{\omega} E_{{\boldsymbol \omega}^*}
\leq\big(\text{max}_{1\leq m \leq N} \text{ diam}_{\omega}E_i\big)\cdot q^{-\frac{k-1}{d}}\\
&=D\cdot q^{1/d}\cdot(q^{-1/d})^k
\leq
D\cdot q^{1/d}\cdot (1/h)^{\frac{d}{\alpha}}\|x_1-x_2\|^{\frac{1}{\alpha}},
\end{split}
\end{equation}
where $D=\max_{1\leq i\leq N} \text{ diam}_{\omega} E_i.$

\begin{figure}[h]
  \includegraphics[width=0.47 \textwidth]{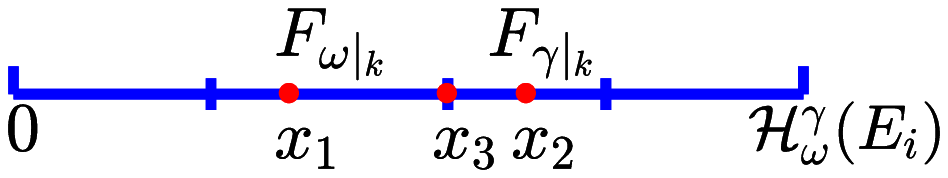}
  \includegraphics[width=0.47 \textwidth]{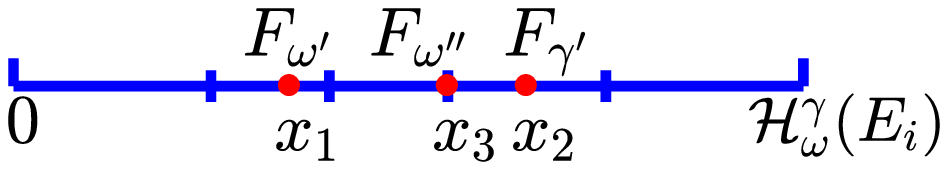}
  \vspace{-1.0 cm}
 \caption{${\bomega}$ and ${\bgamma}$ are adjacent.}\label{adjacent}
\end{figure}

Now, we consider the case that ${\bomega}_k$ and ${\bgamma}_k$ are adjacent. (See Figure \ref{adjacent} (left).)
Let $x_3$ be the intersection of $F_{\bomega|_k}$ and
 $F_{\bgamma|_k}$. Let $k'$ be the smallest integer such that $x_1$ and $x_3$ belong to different cylinders of rank $k'$, say, $x_1\in \rho_i([\bomega'])$ and $x_3\in \rho_i([\bomega''])$ (see Figure \ref{adjacent} (right)), then
 $\|x_1-x_3\|\geq \text{diam}~F_{\bomega''}$
 since $x_3$ is an endpoint.
 Let $y_3=\psi_i(x_3)$. Similar to Case 1,  we have
 $$
 \|y_1-y_3\|_{\omega}\leq
D\cdot q^{1/d}\cdot (1/h)^{\frac{d}{\alpha}}\|x_1-x_3\|^{\frac{1}{\alpha}}.
$$
By the same argument, we have
 $$
\|y_2-y_3\|_{\omega}\leq
D\cdot q^{1/d}\cdot (1/h)^{\frac{d}{\alpha}}\|x_2-x_3\|^{\frac{1}{\alpha}}.
$$
Hence, by the fact $x_3$ is located between $x_1$ and $x_2$,
\begin{equation} \label{Holder-2}
\begin{split}
 \|y_1-y_2\|_{\omega}&\leq \beta \cdot\text{ max}\{\|y_1-y_3\|_{\omega},\|y_3-y_2\|_{\omega}\} \\
&\leq \beta\cdot D\cdot q^{1/d}\cdot (1/h)^{\frac{d}{\alpha}}\cdot\text{max }\{\|x_1-x_3\|^{\frac{1}{\alpha}},\|x_2-x_3\|^{\frac{1}{s}}\}\\
&\leq \beta\cdot D\cdot q^{1/d}\cdot (1/h)^{\frac{d}{\alpha}}\cdot \|x_1-x_2\|^{\frac{1}{\alpha}},
\end{split}
\end{equation}
where the first inequality is from Proposition \ref{Pro-Weak} (iv).

Therefore, \eqref{Holder-1} and \eqref{Holder-2} verify  the $1/\alpha$- H\"older continuity of $\psi_i$.
$\Box$


\bibliographystyle{siam}   
\bibliography{biblio}

%
%
%
%
%
%
%
%
\end{document}